\documentclass[11pt]{amsart}
\usepackage{avant}
\usepackage{fullpage}
\theoremstyle{plain}

\newtheorem{thm}{Theorem}[section]
\newtheorem{lem}[thm]{Lemma}
\newtheorem{prop}[thm]{Proposition}
\newtheorem{cor}[thm]{Corollary}
\newtheorem{remark}[thm]{Remark}

\theoremstyle{definition}
\newtheorem{defn}[thm]{Definition}
\theoremstyle{remark}

\newtheorem{exmp}{Example}[section]

\newcommand\pp{X}
\newcommand\op{\mathcal O_{X}}
\newcommand\mox{\mathcal O_{Y}}
\newcommand\motx{\mathcal O_{\tilde Y}}
\newcommand\opt{\mathcal O_{\tilde X}}
\newcommand\mo{\mathcal O}
\newcommand\lra{\longrightarrow}
\newcommand\ra{\rightarrow}
\newcommand\tpp{\tilde X}
\newcommand\ml{\mathcal L}

\newcommand\OOmega{\overline{\Omega}}
\newcommand\omtp{\OOmega^{1}_{\tilde X}}
\newcommand\omp{\OOmega^{1}_{X}}
\newcommand\CC{\mathbb C}
\newcommand\ZZ{\mathbb Z}
\newcommand\PP{\mathbb P}
\newcommand\FF{\mathbb F}
\newcommand\KK{\mathbb K}
\newcommand\mm{\mathcal M}

\newcommand\mmm{{\mathcal O}_{\pp}(2Y + K_{\pp})}
\newcommand\OO{\mathcal O}

\newcommand\sing{\operatorname{sing}}

\newcommand\Pic{\operatorname{Pic}}
\numberwithin{equation}{section}
\newcommand\md{$\mathfrak{d}$}
\newcommand\mr{$\mathfrak{r}$}
\newcommand{\VVsect}{\vspace*{0ex}}
\begin{document}
\title{  Defect and Hodge numbers of hypersurfaces}
\author{ S\l awomir Rams}

\thanks{Research partially supported by KBN grant no 2P03A 016 25
and the DFG Schwerpunktprogramm "Global methods in complex geometry".}

\email{rams@mi.uni-erlangen.de,  Slawomir.Rams@im.uj.edu.pl}

\begin{abstract}
 We define defect for hypersurfaces with A-D-E singularities
 in
 complex projective normal Cohen-Macaulay fourfolds having some vanishing
properties of Bott-type
and prove
formulae for Hodge numbers of big resolutions
 of such hypersurfaces.
 We
 compute Hodge numbers of Calabi-Yau
manifolds obtained as small resolutions of cuspidal triple sextics and
double octics with higher A$_j$ singularities.
 \end{abstract}
\maketitle

\section{Introduction}
\label{sec:intro}

The starting point of this considerations is the computation of Hodge
numbers of double solids, i.e.
double covers $Y_d$ of the three-dimensional
projective space $\PP_3(\CC)$ branched along a degree-$d$ surface $B_d$.
If the branch divisor $B_d$ is smooth,  then $Y_d$ has always the  Hodge
number $h^{1,1}(Y_d) = 1$.
Let us assume that  the branch locus $B_d$  has $\mu$ ordinary double points as its only singularities.
We call the blow-up of $Y_d$ along $\sing(Y_d)$ the big resolution of $Y_d$ and denote it by $\tilde{Y}_d$.
In \cite{Clemens} Clemens showed the following formula
 for the Hodge number of the big resolution
$$
h^{1,1}(\tilde{Y}_d) = 1 + \mu + \delta \, .
$$
Here the first summand comes from the pull-back of the hyperplane section ${\mathcal O}_{\PP_3}(1)$.
The number $\mu$ is also expected as it counts the number of exceptional divisors in $\tilde{Y}_d$.
The integer $\delta$, called the defect by Clemens, however, is a very subtle invariant of the threefold $Y_d$.
It can  be defined as the number of dependent conditions
imposed on homogenous forms of degree  $(3/2 \cdot d - 4)$ on $\PP_3$ by
the vanishing in the nodes of $B_d$:
$$
\delta := h^{0}({\mathcal O}_{\PP_3}( 3/2 \cdot d - 4) \otimes \mathcal{I}_{\sing(B_d)}) -
\left[h^{0}({\mathcal O}_{\PP_3}( 3/2 \cdot d - 4)) - \mu\right] \, .
$$
Later, the defect was defined for nodal hypersurfaces in $\PP_4$
in \cite{Werner} (see also \cite[Concluding Remarks]{vanstrquint}). Defects of linear systems were also
used to compute Betti numbers of singular hypersurfaces in weighted projective spaces
(see \cite[\S~6.4]{Dimca2}).

Cynk \cite{Cynk} gave another proof of Clemens' formula and generalized
it to ample
three-dimensional hypersurfaces $Y$ with ordinary double points in a smooth projective ambient variety $X$
sharing with the projective space some vanishing properties of Bott-type (\cite[Thm~1]{Cynk2}).
He defined  $V_{Y}$ to be the vector space of the global sections of the line bundle
${\mathcal O}_{X}(2Y) \otimes K_X$ vanishing in the singularities of $Y$. In \cite{Cynk2}
the defect of $Y$ is given by the formula:
$$
\delta_{Y} := \dim V_{Y} -\left[h^{0}({\mathcal O}_{X}(2Y) \otimes K_X)  - \mu \right] \,  .
$$

\noindent By \cite[Thm~1]{Cynk2}, if we assume
$h^{2}(\Omega^{1}_{X}) = h^3(\Omega^{1}_{X}(-Y)) = 0$, then the
Hodge number of $\tilde{Y}$ is
$$
h^{1,1}(\tilde{Y}) = h^{1,1}(X) + \mu + \delta_Y \, .
$$
The double solid $Y_d$ can be embedded as an ample  hypersurface in the weighted projective space
$\PP := \PP(1, 1, 1, 1, d/2)$ in such a way that it does not meet the set $\sing(\PP)$. However, if we resolve
the singularity of $\PP$, the proper transform of $Y_d$ is no longer ample. Thus Clemens' formula cannot be
directly derived from \cite[Thm~1]{Cynk2}.

Here we consider the folowing more general situation:
$\pp$ is a projective normal Cohen-Macaulay fourfold and $Y \subset X$ is a hypersurface
such that $\sing(X) \cap Y = \emptyset$. We assume that $Y$ has only A-D-E singularities, i.e.
for every $P \in \sing(Y)$ there exists  local (analytic) coordinates $x_{1,P}, \ldots, x_{4,P}$ centered at $P$
such that the germ of $Y$ at $P$ is given   by an equation
$$
{\mathfrak n}(x_{1,P}, x_{2,P}, x_{3,P}) + x_{4,P}^2 = 0 \, ,
$$
where ${\mathfrak n}(x_{1,P}, x_{2,P}, x_{3,P})$ is the normal form of the equation of a two-dimensional A-D-E singularity (see the table \eqref{tab-norm}). Let a$_m$ (resp. d$_m$, resp. e$_m$) stand for the number of the singularities of $Y$ of the type A$_m$ (resp.
D$_m$, resp. E$_m$). \\
We define the {\sl big resolution }   $\tilde{\pi} : \tilde{Y}  \ra Y$
as the composition
$\tilde{\pi} =  \sigma_{n} \circ \ldots \circ \sigma_{1}$,
where $ \sigma_{j}:  {\tilde Y}^{j} \ra {\tilde Y}^{j-1}$, for $j = 1, \ldots, n$, is the
blow-up with the center $\sing({\tilde Y}^{j-1})$, ${\tilde Y}^{0} := Y$,
and $ {\tilde Y} = {\tilde Y}^{n}$ is smooth.
The main purpose of this paper is to
 define the defect of a hypersurface with A-D-E singularities and
 obtain a  formula
analogous to \cite[Thm~1]{Cynk2} for the Hodge numbers of the big resolution
$\tilde{Y}$.

The definition of the integer $\mu$ has to be adapted as follows:
$$
\mu_Y := \sum_{m \geq 1} a_m \cdot \lceil m/2 \rceil +  \sum_{m \geq 4} 2 \cdot d_m \cdot \lfloor m/2 \rfloor + 4 \cdot e_6 +
7 \cdot e_7 +  8 \cdot e_8 \, .
$$
Let ${\mathfrak V}_{Y}$ be the space of global sections $H$ of the sheaf ${\mathcal O}_{\pp}(2Y + K_{\pp})$
that vanish in all points $P \in \sing(Y)$ and satisfy the conditions:
\begin{eqnarray*}
\mbox{if } P \mbox{ is an A}_m \mbox{ point, } m \geq 1, &  \hspace*{-13ex} \mbox{ then } \frac{\partial^j H}{ \partial x_{1,P}^j}(P) = 0 \mbox{ for } j \leq \lceil m/2  \rceil -1 \, ,  \\
\mbox{if } P \mbox{ is a D}_m \mbox{ point, } \hspace*{1ex} m \geq 4, &  \hspace*{-2ex} \mbox{ then }  \frac{\partial H}{ \partial x_{2,P}}(P) =
\frac{\partial^j H}{ \partial x_{1,P}^j}(P) =     0   \mbox{ for } j \leq \lfloor m/2 \rfloor - 1 \, , \\
\mbox{if } P \mbox{ is an E}_m \mbox{ point, }  m= 6, & \hspace*{-2ex} 7, 8,  \mbox{ then }
  \frac{\partial H}{ \partial x_{2,P}}(P) =
 \frac{\partial^j H}{ \partial x_{1,P}^j}(P) =     0
  \,  \mbox{ for } j \leq m - 5 \, .
\end{eqnarray*}
The notion of defect has to be adapted in the following way:
$$
\delta_Y :=   \dim({\mathfrak V}_{Y}) - \left[ h^0({\mathcal O}_{\pp}(2Y + K_{\pp})) - \mu_Y \right].
$$
Since we work on a singular ambient variety $X$, we consider  the
Zariski sheaf of germs of $1$-forms
$\OOmega^{1}_{\pp} = j_{*}\Omega^1_{\mbox{\tiny reg}X}$,
where  $j$ stands for the inclusion $\mbox{reg}(X) \ra X$,
 in place of the sheaf of differentials $\Omega^{1}_{\pp}$.
The Bott-type assumptions read as follows:
\begin{description}
\item [{[A1]}] $H^{i}({\mathcal O}_{X}(-Y))= 0$ for $i \leq 3$ \, and \, $H^{j}({\mathcal O}_{X}(-2Y))= 0$
for $j \leq 2$,
\item [{[A2]}] $H^{2}(\OOmega^{1}_{X})=0$,
\item [{[A3]}] $H^{i}(\OOmega^{1}_{X}\otimes \mathcal {\mathcal O}_{X}(-Y))=0, \mbox{ for } i = 1, 2, 3.$
\end{description}
Here we show (Thm~\ref{ade}) that
$$
  h^{1,1}(\tilde Y)=h^{1}(\OOmega^{1}_{X})+\mu_Y +\delta_Y + h^{3}({\mathcal O}_{X}(-2Y)) \, , \\
$$
and, if $h^2(\op) = 0$, then
\begin{align*}
  &h^{1,2}(\tilde Y)=
  h^{0}(\mmm)+h^{4}(\OOmega^{1}_{X})-h^{0}({\mathcal O}_{\pp}(Y + K_{\pp}))-h^{3}(\OOmega^{1}_{X})-
  \\&\rule{1.4cm}{0cm}
    -h^{4}(\OOmega^{1}_{\pp}\otimes{\mathcal O}_{\pp}(-Y))-\mu_Y+\delta_Y.
\end{align*}
In  this way we obtain formulae that
can be applied to  a large class of ambient spaces and to hypersurfaces with higher sigularities.
In particular, Thm~\ref{ade} implies \cite[Cor.~2.32]{Clemens} and \cite[Thm~1]{Cynk2}.

The assumptions [A1], [A2], [A3] are satisfied when $Y$ is an ample hypersurface in a complete simplicial toric fourfold
(Cor.~\ref{projective}).
As an application of Thm~\ref{ade} we derive formulae for the Hodge numbers of various covers of $\PP_3$:
double solids branched along surfaces with Du Val singularities, cyclic n-fold covers
branched along nodal hypersurfaces and triple solids whose branch divisor has A$_2$ singularities.

Sect.~\ref{sect-svb} is devoted to the study of the relation between the Hodge numbers of
a K\"ahler small resolution and the big one.
In the last section of the paper we compute the Hodge numbers of the Calabi-Yau manifolds
obtained as K\"ahler small resolutions of the triple solids branched along the sextics studied in \cite{br1}, \cite{br2}, \cite{labs2}.

\vspace{1ex}
\noindent
{\em Notations and conventions:} All varieties are defined over
the
base-field $\CC$.
 By a divisor we mean a Weil divisor, and "$\sim$" stands for the linear equivalence.
The round-up, resp. the round-down is denoted by $\lceil - \rceil$, resp. $\lfloor - \rfloor$.

\VVsect
\section{Technical preliminaries}
\label{sec:defect}

\newcommand{\Cen}{C}

\noindent
In this section
we assume $\pp$ to be a projective normal Cohen-Macaulay fourfold
and consider a hypersurface $Y \subset \pp$ with {\sl isolated double points} as only singularities.
Moreover, we require
that
\begin{equation} \label{omijaosobliwosci}
\mbox{sing}(\pp) \cap Y = \emptyset \, .
\end{equation}
Here we modify several results from \cite{Cynk2}
in order to apply them in our set-up. Since we are going to work on a
singular (normal) variety $X$, let us recall
that for a Weil divisor $D := \sum_{\Gamma} n_{\Gamma} \,  \Gamma$ on $X$
one defines the sheaf  ${\mathcal O}_X(D)$ by putting
$$
{\mathcal O}_X(D)(U) := \{ f \in \operatorname{Rat}(X) : \operatorname{v}_{\Gamma}(f) +  n_{\Gamma} \geq 0
\mbox{ for every } \Gamma \cap U \neq \emptyset \} \, ,
$$
where  $\operatorname{v}_{\Gamma}(\cdot)$ is the discrete valuation given by the prime divisor $\Gamma \subset X$. Then, the map
$
D \rightarrow  {\mathcal O}_X(D)
$
gives  one-to-one correspondence between the linear equivalence classes of
Weil divisors and isomorphism classes of rank-$1$ reflexive sheaves on $X$ (see \cite[p.~281]{reid1} for the details).

Moreover, if $D_1$ is  Cartier, then one can show that
\begin{equation}
 \label{eq-tens}
{\mathcal O}_X(D_1 + D_2) =  {\mathcal O}_X(D_1) \otimes   {\mathcal O}_X(D_2) \, ,
\end{equation}
(this equality does not hold in general - see \cite[Remark~(5)]{reid1}).

Let $j$ stand for the inclusion $\mbox{reg}(X) \ra X$.
Observe that $\pp$ is endowed with the dualising sheaf $\omega_X = j_{*} \omega_{\mbox{\tiny reg}(X)}$,
 that is reflexive of rank one (see e.g. \cite{reid1}, \cite[Prop.~5.75]{Kollar}), and
we have $\omega_{X} = {\mathcal O}_X (K_X)$. For the definition of the canonical divisor of a normal variety
see e.g. \cite[Def.~0-2-1]{kmm}.

We define
$$
V_{Y} := H^{0}(\mmm\otimes\mathcal I_{\mbox{\tiny sing}(Y)})
$$
to be the vector space of global sections of the
sheaf $\mmm$ on $\pp$
vanishing in all singularities of $Y$.
Let  $\nu$ be the number of points in $\sing(Y)$.
We put
\begin{equation} \label{ry}
r_{Y} =\dim V_{Y} -\left[h^{0}(\mmm)-\nu\right]  .
\end{equation}
Note that the integer
$r_{Y}$ is non-negative because it is  the difference between the actual dimension of the
space  $V_Y$ and the one expected for a hypersurface with
singular points in general position.

Let  $\sigma:\tpp\lra\pp$ be the blow--up of $\pp$ along $\sing(Y)$
and let $\tilde Y$ (resp.  $E$) stand for
the strict transform of $Y$ (resp. the exceptional divisor
of $\sigma$ in $\tpp$).
Since we blow up points in $\mbox{reg}(X)$, we have the equalities:
 \begin{eqnarray}
\sigma_{\ast}(\opt(-E))= {\mathcal I}_{\mbox{\tiny sing}(Y)}, &R^{i}\sigma_{\ast}(\opt(-E))=0, & \mbox{ for }   i>0,  \label{directimageofmE}  \\
\sigma_{\ast}(\opt(kE))= \op, \hspace{4ex} &R^{i}\sigma_{\ast}(\opt(kE))=0,
& \mbox{ for }  k = 0, 1, 2, 3 \mbox{ and }  i>0, \label{directimageofkE} \\
\sigma_{\ast}(\opt(4E))= \op, \hspace{4ex}
 &R^{j}\sigma_{\ast}(\opt(4E))=0,\hspace{0ex}  &
\mbox{ for } j = 1,2, \label{directimageoftE} \\
 R^{3}\sigma_{\ast}(\opt(4E))=\CC^{\nu},  \hspace{2ex} & &  \nonumber
\end{eqnarray}
where $\CC^{\nu}$ stands for the skyscraper sheaf with stalk $\CC$ over the centers of blow-up. \\
Indeed, observe that
the first equality in \eqref{directimageofmE} and the
equalities  \eqref{directimageofkE}  for $k = 0$ are
 basic properties of
blow-up. To prove the second assertion of   \eqref{directimageofmE} apply the direct image $\sigma_{\ast}$ to the exact sequence
\begin{equation} \label{exc-div}
    0\lra\opt(-E)\lra\opt\lra
 {\mathcal O}_{E}  \lra0 \, .
\end{equation}
To show  \eqref{directimageofkE} for $k = 1,2,3$ and  \eqref{directimageoftE}
 consider the tensor
product of \eqref{exc-div} with the locally free sheaf $\opt(kE)$
(see \cite[the proof of Lemma~1]{Cynk2} for details).

One can see that
if ${\mathcal F}$ is a coherent sheaf on $\tpp$ and ${\mathcal E}$ is a coherent sheaf on $\pp$ such that
 ${\mathcal E}|_{\mbox{\tiny reg}(X)}$ is locally free, then
\begin{equation} \label{uogolnioneprojectionformula}
R^{i}\sigma_{*}({\mathcal F} \otimes_{{\mathcal O}_{\tpp}}  \sigma^{*}{\mathcal E}) =
 R^{i}\sigma_{*}{\mathcal F} \otimes_{{\mathcal O}_{\pp}} {\mathcal E} \,
\mbox{ for } i \geq 0.
\end{equation}

\noindent
Indeed, the blow-up $\sigma$, when restricted to an appropriate neighbourhood
of singular points of $\pp$ and $\tpp$, is  an isomorphism, so the stalks
of both sheaves vanish for a point in such a neighbourhood and $i \geq 1$.
To complete the proof of  \eqref{uogolnioneprojectionformula}
for $i \geq 1$
apply the projection formula \cite[III.~Ex.~8.3]{hart}
on $\mbox{reg}(X)$.
Finally, by the assumption \eqref{omijaosobliwosci},
one can literally repeat the proof of
the projection formula
(see e.g. \cite[Lemma~2.29]{Iitaka}) to show that
the formula in question holds when $i = 0$. \\
In particular, we get
 \begin{equation} \label{slabaproj}
\sigma_{*}(\sigma^{*}{\mathcal O}_{\pp}(K_{\pp}) \otimes_{{\mathcal O}_{\tilde X}}
{\mathcal O}_{\tilde X}(-E))) = {\mathcal O}_{\pp}(K_{\pp}) \otimes_{{\mathcal O}_{X}} {\mathcal I}_{\mbox{\tiny sing}(Y)} \, .
\end{equation}
Let  $D$ be a  Weil divisor on $X$.
Since the centers of the blow-up $\sigma$ are smooth on $\pp$,
the pull-back $\sigma^{*}D$ is well-defined.
Moreover,  if we choose a divisor $K_{\pp}$ (resp. $K_{\tpp}$) in the canonical class of $\pp$
(resp. $\tpp$), then we have
\begin{equation} \label{canonical}
K_{\tpp} \sim \sigma^{*}K_{\pp} + 3E .
\end{equation}

In the rest of this section  we put
$$
 \mm:=\mo_{X}(Y) \,, \, \, \, \ml := \opt(\tilde Y) \, .
$$

\begin{lem}\label{formulanadefekt} Let
$\tilde{Y}$
be the  proper transform of
$Y$
under
the blow-up
with center $\sing(Y)$.
If
$$
h^{i}(\mm^{-j})= 0 \mbox{ for } i \leq 2 , j = 1,2, \mbox{ and } h^{3}(\mm^{-1})= 0,
$$
then \\
a)  $h^{i}(\ml^{-j})= 0 \mbox{ for } i \leq 2 , j = 1,2,$ and $h^{3}(\ml^{-1})= 0$, \\
b)  $h^{i}({\mathcal O}_{\tilde{Y}}) =   h^{i}(\op), \mbox{ for } i \leq 2$, \\
c)  $h^{0}(\motx\otimes\ml^{-1})=h^{1}(\motx\otimes\ml^{-1})=0, $ \\
d)  $h^{4}(\ml^{-2}) =  \dim V_{Y}$, \\
e)  $h^{3}(\ml^{-2}) = h^{3}(\mm^{-2}) +  r_{Y}.$
\end{lem}
\begin{proof}
By~(\ref{omijaosobliwosci}) the sheaf  $\ml$ is locally free, so
we can apply the projection formula and \eqref{directimageofkE} to
the bundle  $\ml^{-1}=\sigma^{*}\mm^{-1}\otimes \opt(2E)$. In this
way we obtain  the equalities
  \begin{eqnarray} \label{directimageofl}
 \sigma_{\ast}(\ml^{-1})=\mm^{-1},&  R^{i}\sigma_{\ast}(\ml^{-1})=0 \mbox{ for } i>0. &
 \end{eqnarray}
Similarly, \eqref{directimageoftE} yields
\begin{eqnarray} \label{directimageofl2}
 \sigma_{\ast}(\ml^{-2})=\mm^{-2},&  R^{i}\sigma_{\ast}(\ml^{-2})=0,  \mbox{ for } i=1,2,
    &\mbox{and } R^{3}\sigma_{\ast}(\ml^{-2})=\CC^{\nu}.
\end{eqnarray}
Moreover, by \eqref{directimageofmE} and \eqref{uogolnioneprojectionformula} we have
\begin{eqnarray} \label{directimageoftwist}
 \sigma_{\ast}(\sigma^{*}\mm^{2} \otimes \sigma^{*}{\mathcal O}_{\pp}(K_{\pp}) \otimes
{\mathcal O}_{\tpp}(-E)) & = & \mm^{2}  \otimes {\mathcal O}_{\pp}(K_{\pp}) \otimes \mathcal I_{\mbox{\tiny sing}(Y)},  \label{directimageoftwist1} \\
  R^{i}\sigma_{\ast}(\sigma^{*}\mm^{2} \otimes \sigma^{*}{\mathcal O}_{\pp}(K_{\pp}) \otimes
{\mathcal O}_{\tpp}(-E)) & = & 0,  \mbox{ for } i \geq 1.  \label{directimageoftwist2}
\end{eqnarray}

The Leray spectral sequence for the map $\sigma$ and the sheaf $\ml^{-1}$, combined with \eqref{directimageofl},
yields for all $i$
\begin{equation} \label{vanl1}
 h^{i}(\ml^{-1})=h^{i}(\mm^{-1}) \, .
\end{equation}
Let us consider the  Leray spectral sequence for
$\sigma$ and $\ml^{-2}$.
Since $E^{p,q}_2 = H^p(X, R^q\sigma_{\ast} \ml^{-2})$, we have
$
E^{0,3}_2 = \CC^{\nu}
$
and, by \eqref{directimageofl2},  the $E^{p,q}_2$-terms vanish  for $q > 0$ and $q \neq 3$.
Therefore, (see e.g. \cite[Ex.~1.D]{mccleary}) we have
 the equalities
\begin{equation} \label{vanl2}
h^{i}(\ml^{-2})=h^{i}(\mm^{-2})  \mbox{ for } i \leq 2
\end{equation}
and the  exact sequence
\begin{equation} \label{leraydlalm}
0 \lra H^{3}(\mm^{-2})\lra H^{3}(\ml^{-2})\lra \CC^{\nu} \lra
H^{4}(\mm^{-2})\lra H^{4}(\ml^{-2})\lra0 \, .
\end{equation}
Finally,  we have
\begin{eqnarray} \label{eq-6}
h^{0}(\ml^{2} \otimes {\mathcal O}_{\tpp}(K_{\tpp})) & = &   h^{0}(\sigma^{*}\mm^{2} \otimes  {\mathcal O}_{\tpp}(-4E) \otimes
{\mathcal O}_{\tpp}(K_{\tpp})) \nonumber \\
&  \stackrel{\eqref{canonical}}{=} &
h^{0}(\sigma^{*}\mm^{2} \otimes \sigma^{*}{\mathcal O}_{\pp}(K_{\pp}) \otimes
{\mathcal O}_{\tpp}(-E)) \\
&  = &
 h^{0}(\mm^{2} \otimes {\mathcal O}_{\pp}(K_{\pp})  \otimes  {\mathcal I}_{\mbox{\tiny sing}(Y)} ),  \nonumber
\end{eqnarray}
where the last equality results from the Leray spectral sequence, \eqref{directimageoftwist1} and   \eqref{directimageoftwist2}.

Now we are in position to prove the claims a) - e). \\
{\sl a)} results immediately from \eqref{vanl1}, \eqref{vanl2}. \\
{\sl b)} Obviously $h^{i}(\opt) =   h^{i}(\op)$ for all $i$. For  $i \leq 2$,
the equality $h^{i}({\mathcal O}_{\tilde{Y}}) =   h^{i}(\opt)$
results from a) and the long cohomology sequence associated to the short exact sequence
$$
  0\lra\ml^{-1}\lra\mo_{\tpp}\lra\motx\lra0.
$$
{\sl c)} results from a) and the long cohomology sequence associated to the  sequence
$$
  0\lra\ml^{-2}\lra\ml^{-1}\lra\motx\otimes\ml^{-1}\lra0 \, .
$$
{\sl d)} The sheaf $\ml^{-2}$ is locally free and $\tpp$ is
Cohen-Macaulay, so we can use Serre duality
(\cite[Cor.~7.7]{hart}) to prove that
\begin{equation} \label{eq-1-pico}
h^{4}(\ml^{-2}) = h^{0}(\ml^{2} \otimes {\mathcal O}_{\tpp}(K_{\tpp}))
  \stackrel{\eqref{eq-6}}{=}  \dim V_{Y}  \, .
\end{equation}
{\sl e)} By the same argument we have
\begin{equation}  \label{eq-l1-2p}
 h^{0}(\mm^{2} \otimes {\mathcal O}_{\pp}(K_{\pp}))  =  h^{4}(\mm^{-2}) .
\end{equation}
From the exact sequence
$$
  0\lra\mm^{-2}\lra\mm^{-1}\lra\motx\otimes\mm^{-1}\lra0 \, ,
$$
we obtain
\begin{eqnarray}
 h^{2}(\mox\otimes\mm^{-1}) & = & h^{3}(\mm^{-2})  \label{eq-l1-2}  \, ,  \\
 h^{3}(\mox\otimes\mm^{-1}) & = & h^{4}(\mm^{-2}) -  h^{4}(\mm^{-1}) \label{eq-l1-2pp}  \, .
\end{eqnarray}
Finally, we have the equalities
\begin{eqnarray*}
  h^{3}(\mm^{-2}) +  r_{Y}
&  \stackrel{\eqref{eq-1-pico}}{=} &  h^{3}(\mm^{-2}) +  h^{4}(\ml^{-2}) +
                       \nu -  h^{0}(\mm^{2} \otimes {\mathcal O}_{\pp}(K_{\pp})) \\
  & \stackrel{\eqref{eq-l1-2p}}{=}  &  h^{3}(\mm^{-2}) +  h^{4}(\ml^{-2}) +  \nu -   h^{4}(\mm^{-2}) \stackrel{\eqref{leraydlalm}}{=}   h^{3}(\ml^{-2}) \, .
 \end{eqnarray*}
\end{proof}
\noindent
It should be pointed out that if $X$ has isolated rational singularities and $\mm$ is ample,
then
the assumptions of Lemma~\ref{formulanadefekt} are satisfied by  \cite[VII.Thm~7.80]{ss}.


We intend to apply \cite[Bott's Vanishing,~p.~130]{oda}
(see also \cite[Thm~2.3.2]{dol})
in Sect.~\ref{sec:brts}, so in the remainder of this section we
focus our attention on the Zariski sheaf of germs of $1$-forms $
\OOmega^{1}_{\pp} := j_{*}\Omega^1_{\mbox{\tiny reg}X},$
 where $j:\mbox{reg}(X) \rightarrow X$ stands for the inclusion.
Recall that the sheaf $\OOmega^{1}_{\pp}$ does not have to be
locally free.

Let  $C_l$, where $l = 1, \ldots, \nu$,
stand for a  component of the exceptional divisor $E$.
We have the following exact sequence
\begin{equation} \label{ciagdoklrr}
    0\lra\sigma^{\ast}\omp\lra\omtp\lra
 \bigoplus\limits_{l = 1}^{\nu}\Omega^{1}_{C_{l}}  \lra0 \, ,
\end{equation}
Indeed, away from the singularities of $X$ the shaves $\OOmega^{1}_{X}$ and $\Omega^{1}_{X}$
coincide, whereas around the singular locus $\sing(X)$ the blow-up is an isomorphism. \\


\begin{lem} \label{dirimpullbom} We have
 $ \sigma_{\ast} ( \sigma^{\ast}\omp ) \cong \omp $
and  $R^{i} \sigma_{\ast} (\sigma^{\ast} \OOmega^{1}_{\pp}) = 0$ for $i > 0$.
\end{lem}
\begin{proof}
Apply  \eqref{uogolnioneprojectionformula} and \eqref{directimageofkE} for $k = 0$.
\end{proof}

Now, we are in position to prove

\begin{lem}\label{hiomegatenslminone} a) If $h^2(\OOmega^{1}_{\pp}) = 0$, then  $h^{1}\omtp={h^{1}\omp+\nu}$ and
 $h^{i}\omtp={h^{i}\omp}$ for $ i \not=1$. \\
b) The equality $h^{i}(\omtp\otimes\ml^{-1})=h^{i}(\omp\otimes\mm^{-1})$ holds for every $i$.
\end{lem}
\begin{proof} a) (c.f. \cite[Lemma~2]{Cynk2})
Recall that all centers of the blow-up $\sigma$ are smooth on $X$,
so $C_l \simeq  \PP_3$ for $l = 1, \ldots, \nu$.
Apply the direct image  $\sigma_{\ast}$ to the  Euler sequence
to see that
$R^{1}\sigma_{\ast}\Omega^{1}_{C_{l}}$
is the skyscraper sheaf with stalk $\CC$ at each center of the blow-up, whereas
$R^{i}\sigma_{\ast}\Omega^{1}_{C_{l}} = 0 \mbox{ for } i \neq 1.$
Consider the direct image  $\sigma_{\ast}$ of the exact
sequence \eqref{ciagdoklrr}. Lemma~\ref{dirimpullbom} yields
\begin{eqnarray} \label{eq-spec}
\sigma_{\ast}\OOmega^{1}_{\tpp} = \OOmega^{1}_{\pp}, &
R^{1}\OOmega^{1}_{\tpp} = \CC^{\nu},  &
R^{i}\OOmega^{1}_{\tpp} = 0 \mbox{ for } i \geq 2.
\end{eqnarray}
Now we consider the Leray spectral sequence
$
E^{p,q}_2 = H^{p}(\pp,R^{q}\OOmega^{1}_{\tpp}) \Rightarrow
H^{p+q}(\tpp, \OOmega^{1}_{\tpp}) \, .
$
By \eqref{eq-spec} we have the vanishing
$$
E_{2}^{2,0} = H^2(\sigma_{\ast}\OOmega^{1}_{\tpp}) = H^2(\OOmega^{1}_{\pp}) = 0
$$
The latter implies that the differential
$d_2 \, : \, E_{2}^{0,1} \rightarrow E_2^{2,0}$
is the zero map and we can  compute all $E_{\infty}^{p,q}$-terms. \\
b) In view of the Leray spectral sequence $ E^{p,q}_2 =
H^{p}(\pp,R^{q}(\OOmega^{1}_{\tpp}\otimes\ml^{-1})) \Rightarrow
H^{p+q}(\tpp, \OOmega^{1}_{\tpp}\otimes\ml^{-1}) \, , $ it
suffices to show that
\begin{equation} \label{eq-exc-prep}
 \sigma_{\ast} (\omtp\otimes\ml^{-1}) = \omp\otimes\mm^{-1} \mbox{ and }
R^{i} \sigma_{\ast} (\omtp\otimes\ml^{-1}) = 0 \mbox{ for }i > 0.
\end{equation}
Tensoring \eqref{ciagdoklrr} with $\ml^{-1}$ we get
\begin{equation} \label{eq-exc-dm}
0\lra\sigma^{\ast}\omp\otimes \sigma^{\ast}\mm^{-1}\otimes\opt(2E)\lra\omtp\otimes\ml^{-1}\lra
\bigoplus_{l = 1}^{\nu}\Omega^{1}_{C_{l}}(-2)\lra0.
\end{equation}
For $i \geq 0$, we have the vanishing
$
R^{i} \sigma_{\ast} \Omega^{1}_{C_{l}}(-2)=0 \, ,
$
so applying the
direct image functor to \eqref{eq-exc-dm} yields the exact sequences
\begin{equation} \label{eq-eq-dorz}
0\lra R^{i} \sigma_{\ast} ( \sigma^{\ast}\omp\otimes  \sigma^{\ast}\mm^{-1}
\otimes\opt(2E))\lra R^{i} \sigma_{\ast}(\omtp\otimes\ml^{-1}) \lra0 .
\end{equation}

Now \eqref{uogolnioneprojectionformula}, \eqref{directimageofkE} and Lemma~\ref{dirimpullbom} imply
$
\sigma_{\ast} ( \sigma^{\ast}(\omp)\otimes\opt(2E)) \cong \omp \, .
$
From the projection formula we obtain
$$
 \sigma_{\ast} ( \sigma^{\ast}(\omp)\otimes \sigma^{\ast}(\mm^{-1})\otimes\opt(2E)) \cong
\mm^{-1} \otimes  \sigma_{\ast} ( \sigma^{\ast}(\omp)\otimes\opt(2E)) \, ,
$$
which completes the proof of the first claim in \eqref{eq-exc-prep}.

Let $i > 0$.
To prove the second claim of \eqref{eq-exc-prep} apply
 \eqref{uogolnioneprojectionformula} and
\eqref{directimageofkE} for $k=2$ to show that
$$
R^i \sigma_{\ast}( \sigma^{\ast}\omp\otimes  \sigma^{\ast}\mm^{-1}
\otimes\opt(2E)) = 0.
$$
and use \eqref{eq-eq-dorz}.
\end{proof}

\begin{lem}\label{hiomegatensoy}
If $h^2(\OOmega^{1}_{\pp}) = 0$ and  $h^{i}(\omp\otimes\mm^{-1}) = 0$ for $i = 1,2,3$, then
  \[
  \begin{array}{l}
    h^{0}(\omtp\otimes \motx)=h^{0}(\OOmega^{1}_{X})-h^{0}(\omp\otimes\mm^{-1}),   \\
    h^{1}(\omtp\otimes \motx)=h^{1}(\OOmega^{1}_{X})+\nu,\\
    h^{2}(\omtp\otimes \motx)=0, \\
    h^{3}(\omtp\otimes \motx)=h^{3}(\OOmega^{1}_{X})
+h^{4}(\OOmega^{1}_{X}\otimes \mm^{-1})-h^{4}(\OOmega^{1}_{X}) \, .
  \end{array}\]
\end{lem}
\begin{proof}
Observe that  the following short
sequence
\[0\lra\omtp\otimes \ml^{-1}\lra\omtp\lra\omtp\otimes\motx\lra0\]
is exact. Apply Lemma~\ref{hiomegatenslminone}
to the associated long cohomology sequence.
\end{proof}

\VVsect
\section{Defect for A-D-E singularities}
\label{sec:brvsr}
In this section we define the {\sl defect} for threefold hypersurfaces with
A-D-E singularities and give a  formula for this
integer. We work with semiquasihomogenous equations of A-D-E germs, instead of their normal forms, because the analytic
coordinates in which a germ is given by the former are easier to find, while the formula for defect remains the same.

Let $X$ be a projective normal Cohen-Macaulay fourfold and let $Y \subset X$ be a hypersurface
with A-D-E singularities
such that $\mbox{sing}(\pp) \cap Y = \emptyset$.
In particular the singularities of $Y$ are absolutely isolated, i.e. $\sing(Y)$ can be resolved by blowing up (closed) points
(\cite[p.137]{durfee}).

We define the {\sl big resolution }   $\tilde{\pi} : \tilde{Y}  \ra Y$
as the composition:
\begin{equation} \label{eq-big}
\tilde{\pi} =  \sigma_{n} \circ \ldots \circ \sigma_{1} \, : \,  \tilde{Y} \ra Y = : \tilde{Y}^{0} \, ,
\end{equation}
where $ \sigma_{j}:  {\tilde Y}^{j} \ra {\tilde Y}^{j-1}$, for $j = 1, \ldots, n$, is the blow-up with the center
$\sing({\tilde Y}^{j-1}) \neq \emptyset$
and $ {\tilde Y} = {\tilde Y}^{n}$ is smooth.
By abuse of notation we use the same symbol to denote the composition of the  blow-ups of $X$
with the same centers
$
\tilde{\pi}: \tilde{X} \ra X \, .
$ \\
We put $\tilde{Y} = Y$, $n = 0$, when $Y$ is smooth.
The number of singular and infinitely
near singular points of $Y$ is denoted by
$$
\mu_{Y} := \sum_{j=0}^{n-1} \# \sing({\tilde Y}^{j}) \, .
$$
Since  the singularities of ${\tilde Y}^{j}$ are isolated for
every $j$, we can formulate the following definition
\begin{defn} \label{defect}
We define the \emph{defect of $Y$} as the non--negative number
  \[\delta_{Y} = r_{Y} + \ldots + r_{\tilde{Y}^{n-1}} \, . \]
\end{defn}
\noindent Observe that if $Y$ has only ordinary double points and
$X$ is smooth, then this definition coincides with
\cite[Def.~1]{Cynk2}.

Let $j = 1, \ldots, n$ and
let $E_j$ stand for the exceptional divisor of the blow-up
$\sigma_{j}:  {\tilde X}^{j} \ra {\tilde X}^{j-1}$.
For a divisor $H$ on $\tilde{X}^{j-1}$, we define the (Weil) divisor
$$
{\mathfrak s_j}(H) := \sigma_{j}^{*}H - E_j.
$$

To simplify our notation, in the following lemma  we use the same letter to denote a section
of $\mmm$
and the divisor it defines.
\begin{lem} \label{uproszdef}
Let $\sing(Y)$ consist of A-D-E singularities and let $\tilde{\pi} =  \sigma_{n} \circ \ldots \circ \sigma_{1}$
be the big resolution of $Y$. Then
\begin{eqnarray*}
\delta_{Y} &=& h^{0}({\mathcal O}_{\tilde{X}}(2 \tilde{Y} + K_{\tilde{X}})) - h^{0}(\mmm) + \mu_{Y} \, ,  \\
 h^{0}({\mathcal O}_{\tilde{X}}(2 \tilde{Y} + K_{\tilde{X}}))   &=& \dim \{ H \in V_{Y}  \, : \,
\operatorname{sing}(\tilde{Y}^{j}) \subset  \operatorname{supp}(({\mathfrak s_j} \circ \ldots \circ {\mathfrak s_1})(H))
\mbox{ for } j  \leq n-1\}.
\end{eqnarray*}
\end{lem}
\begin{proof}
For $j = 1, \ldots, n$ the singularities of $\tilde{Y}^{j-1}$ are double points. In this case the equality
\eqref{eq-6} (see also \eqref{eq-tens}) reads
\begin{equation} \label{eq-cof}
h^{0}({\mathcal O}_{\tilde{X}^{j}}(2 \tilde{Y}^{j} + K_{\tilde{X}^{j}}))
=
h^{0}({\mathcal O}_{\tilde{X}^{j-1}}(2 \tilde{Y}^{j-1} + K_{\tilde{X}^{j-1}})
\otimes\mathcal I_{\mbox{\tiny sing}( \tilde{Y}^{j-1})}) \, .
\end{equation}
Now the formula for $\delta_Y$ follows directly from Def.~\ref{defect}.

The second equality is obvious for $n = 1$, so we can assume $n \geq 2$.
Fix $j = 1, \ldots, n$ and
observe that  the natural linear map
$$
V_{\tilde{Y}^{j-1}} \lra
H^{0}({\mathcal O}_{\tilde{X}^{j}}(\sigma_{j}^{*}(2 \tilde{Y}^{j-1}
+ K_{\tilde{X}^{j-1}}) - E_{j})) =
H^{0}({\mathcal O}_{\tilde{X}^{j}}(2 \tilde{Y}^{j} + K_{\tilde{X}^{j}})) \,
$$
that maps a section which defines a divisor $D \in |2 \tilde{Y}^{j-1} + K_{\tilde{X}^{j-1}}|$
satisfying the condition
$$
\mbox{sing}(\tilde{Y}^{j-1}) \subset  \mbox{supp}(D)
$$
  to its lift, which is a section that gives the divisor ${\mathfrak s_j}(D)$, is an isomorphism. We obtain
the equalities
\begin{eqnarray*}
 h^{0}({\mathcal O}_{\tilde{X}}(2 \tilde{Y} + K_{\tilde{X}})) &\stackrel{\eqref{eq-cof}}{=}&   \dim V_{\tilde{Y}^{n-1}} \\
 &=& \dim \{ H \in  V_{\tilde{Y}^{n-2}} \, : \,
\mbox{sing}(\tilde{Y}^{n-1}) \subset  \mbox{supp}({\mathfrak s_{n-1}}(H)) \} \, = \ldots \\
  &=&   \dim \{ H \in V_{Y}  :
\mbox{sing}(\tilde{Y}^{j}) \subset  \mbox{supp}(({\mathfrak s_j} \circ \ldots \circ {\mathfrak s_1})(H))
\mbox{ for }
j  \leq n-1\}
\end{eqnarray*}
\end{proof}
\noindent
Remark: In the proof of Lemma~\ref{uproszdef} we used the apparently weaker assumption that
 all singularities and infinitely near singularities of $Y$ are isolated double points.
By \cite[Thm~1]{treger} if $P \in Y$ is
an absolutely isolated double point, then it is  an A-D-E singularity.
\vspace*{1ex}

For an A-D-E point $P \in \sing(Y)$ we choose (analytic) coordinates $x_{1,P}, \ldots, x_{4,P}$ centered at $P$
such that the germ of $Y$ at $P$ is given by the
equation
\begin{equation} \label{eq-sqh}
{\mathfrak n}(x_{1,P}, x_{2,P}, x_{3,P}) + x_{4,P}^2 + F(x_{1,P}, x_{2,P}, x_{3,P}, x_{4,P}) = 0 \, ,
\end{equation}
where ${\mathfrak n}(x_{1,P}, x_{2,P}, x_{3,P})$ is the normal form of the equation of a two-dimensional A-D-E singularity and
$F(x_{1,P}, x_{2,P}, x_{3,P}, x_{4,P})$ is
 a polynomial of order strictly greater than $1$ with respect to the weights
$\operatorname{w}_{\mathfrak n}(x_{1,P})$, $\operatorname{w}_{\mathfrak n}(x_{2,P})$, $\operatorname{w}_{\mathfrak n}(x_{3,P})$,
$\operatorname{w}_{\mathfrak n}(x_{4,P})$
 given in the table below:
\vspace*{-0.2ex}
\begin{equation} \label{tab-norm}
\renewcommand\arraystretch{1.5}
\begin{tabular}{|c|c|c|}
\hline
          & ${\mathfrak n}(x_{1}, x_{2}, x_{3})$ &  ($\operatorname{w}_{\mathfrak n}(x_{1}), \ldots , \operatorname{w}_{\mathfrak n}(x_{4})) $  \\ \hline
A$_m$,  $m \geq 1$  &  \hspace{1ex}   $x_{1}^{m+1} + x_{2}^2 + x_{3}^2$ \hspace{1ex} & $(\frac{1}{m+1}, \frac{1}{2}, \frac{1}{2}, \frac{1}{2})$ \\ \hline
D$_m$,  $m \geq 4$ &   \hspace{1ex}  $x_{1} \cdot (x_{2}^2 + x_{1}^{m-2})  + x_{3}^2$  \hspace{1ex} &
 \hspace{1ex} $(\frac{1}{m-1}, \frac{m-2}{2(m-1)}, \frac{1}{2}, \frac{1}{2})$  \hspace{1ex} \\ \hline
E$_6$ &  $x_{1}^4 + x_{2}^3 + x_{3}^2$ & $(\frac{1}{4},  \frac{1}{3}, \frac{1}{2}, \frac{1}{2})$  \\ \hline
E$_7$ &  $x_{1}^3 \cdot x_{2} + x_{2}^3 + x_{3}^2$   & $(\frac{2}{9}, \frac{1}{3}, \frac{1}{2}, \frac{1}{2})$  \\ \hline
E$_8$ &  $x_{1}^5 + x_{2}^3 + x_{3}^2$ &  $(\frac{1}{5}, \frac{1}{3}, \frac{1}{2}, \frac{1}{2})$    \\ \hline
\end{tabular}
\renewcommand\arraystretch{1}
\end{equation}

\noindent Recall that, according to \cite[Char.~C~9]{durfee},
every germ given by an equation of the type \eqref{eq-sqh} is
A-D-E, so we can assume  $F = 0$ to check that the number of
singular points (different from  $P$) which are infinitely near
$P$ is as follows \vspace*{-0.2ex}
\begin{equation} \label{tab-infp}
\mbox{
 \renewcommand\arraystretch{1.5}
\begin{tabular}{|c|c|c|c|c|}
\hline
 \hspace{1ex} A$_m$, $m \geq 1$ \hspace{1ex} &  \hspace{1ex}  D$_m$, $m \geq 4$ \hspace{1ex}
& \hspace{2ex} E$_6$ \hspace{2ex} &  \hspace{2ex} E$_7$  \hspace{2ex}
&  \hspace{2ex} E$_8$  \hspace{2ex} \\ \hline
   $\lceil m/2 \rceil - 1$  &     $2 \cdot \lfloor m/2 \rfloor - 1$  &   $3$    &   $6$    &   $7$     \\ \hline
\end{tabular} }
 \renewcommand\arraystretch{1}
\end{equation}
\vskip2mm

Now, we are in position to state a lemma that, combined with Lemma~\ref{uproszdef},
gives a method of computing $\delta_Y$ without studying the configuration of
infinitely near singular points.

\begin{lem} \label{lem-defekt}
Let $\sing(Y)$ consist of  A-D-E points and let  $H \in  H^{0}(\mmm\otimes\mathcal I_{\mbox{\tiny sing}(Y)})$.
Then, the following conditions are equivalent:
\vskip2mm
\noindent
{\rm [\,I\,]}
The inclusion  $\operatorname{sing}(\tilde{Y}^{j}) \subset  \operatorname{supp}(({\mathfrak s_j} \circ \ldots \circ {\mathfrak s_1})(H))$ holds for  $j  \geq 1$ .

\vskip2mm
\noindent
{\rm [II]}
 Every $P \in \sing(Y)$ satisfies one of the following:
\begin{eqnarray}
& & P \mbox{ is an A}_m \mbox{ point, with }m \geq 1 \mbox{ and }
\frac{\partial^j H}{ \partial x_{1,P}^j}(P) = 0 \mbox{ for } j \leq \lceil m/2  \rceil -1 \, ,   \label{eq-cond-am}   \\
& & P \mbox{ is a D}_m \mbox{ singularity,  where }m \geq 4 \mbox{ and }
 \frac{\partial H}{ \partial x_{2,P}}(P) =
\frac{\partial^j H}{ \partial x_{1,P}^j}(P) =     0   \mbox{ for } j \leq \lfloor m/2 \rfloor - 1 \, ,  \label{eq-cond-dm} \\
& & P \mbox{ is an E}_m \mbox{ point,  where }m = 6, 7, 8 \mbox{ and }
 \frac{\partial H}{ \partial x_{2,P}}(P) =
\frac{\partial^j H}{ \partial x_{1,P}^j}(P) =     0
 \,  \mbox{ for } j \leq m - 5 \, , \label{eq-cond-em}
\end{eqnarray}
where  $x_{1,P}, \ldots, x_{4,P}$ are analytic local coordinates centered at the point $P$
such that the hypersurface $Y$ is given near $P$ by the semiquasihomogenous equation \eqref{eq-sqh}.
\end{lem}
\begin{proof} We choose a point $P$ and fix the coordinates $x_{1,P}, \ldots, x_{4,P}$.
 Let
 $\tilde{\pi}  =  \sigma_{n} \circ \ldots \circ \sigma_{1}$ be the big resolution of $Y$.
We claim that for a point $P$ of the type A$_m$ (resp. D$_m$, resp. E$_m$),
 the condition \eqref{eq-cond-am}
(resp. \eqref{eq-cond-dm}, resp. \eqref{eq-cond-em})
is fulfilled iff for all  $j$ such that
$\operatorname{sing}(\tilde{Y}^{j}) \neq \emptyset$ we have the inclusion
\begin{equation} \label{eq-inf-eq}
 \operatorname{sing}(\tilde{Y}^{j}) \cap  (\sigma_{j} \circ \ldots \circ \sigma_{1})^{-1}(P)
 \subset  \operatorname{supp}(({\mathfrak s_j} \circ \ldots \circ {\mathfrak s_1})(H)) \, .
\end{equation}

Suppose that $P$ is a D$_m$ singularity and put $x_j := x_{j,P}$.
To simplify our notation we denote the local equation of the divisor $H$ near the point $P$ by $H$.
Let
$$
F(x_{1}, \ldots, x_4) = \sum f_{j_1, j_2, j_3, j_4} \cdot x_1^{j_1} \ldots x_4^{j_4} \mbox{ and }
H(x_{1}, \ldots, x_4) = \sum h_{j_1, j_2, j_3, j_4} \cdot x_1^{j_1} \ldots x_4^{j_4}  \, .
$$
We consider three cases: \\
{\bf D$_4$}: In this case,  $\tilde{Y}^{1} \subset  X \times \PP_3$ has three A$_1$ points as its only singularities.
Let $(y_1:y_2:y_3:y_4)$ stand for homogenous coordinates on  $\PP_3$. Then in the affine set $y_2 = 1$
the blow-up  $\tilde{Y}^{1}$ is given by the equation
$$
x_2 \cdot y_1 \cdot (1 + y_1^2) + y_3^2 + y_4^2 +
\sum f_{j_1, j_2, j_3, j_4} \cdot y_1^{j_1} x_2^{j_1 + \ldots + j_4 - 2} y_3^{j_3} y_4^{j_4} \, .
$$
Observe that if
$$
1/3 \cdot j_1 + 1/3 \cdot j_2 + 1/2 \cdot j_3 + 1/2 \cdot j_4 > 1,
$$
then the monomial
$y_1^{j_1} \cdot x_2^{j_1 + \ldots + j_4 - 2} \cdot y_3^{j_3} \cdot y_4^{j_4}$ is singular along the set
$x_2 =  y_3 =  y_4 = 0$. Consequently,  for every choice of
$F(x_1, \ldots, x_4)$ in the semiquasihomogenous equation \eqref{eq-sqh},
the variety   $\tilde{Y}^{1}$ is singular in the points
$$
P_{1,1} := (P,(0:1:0:0)), P_{1,2} := (P,(i:1:0:0)), P_{1,3} := (P,(-i:1:0:0))
\in  X \times \PP_3.
$$
For $y_2 = 1$,  the divisor  ${\mathfrak s_1}(H)$ is given by
$$
 \sum h_{j_1, j_2, j_3, j_4} \cdot y_1^{j_1}   x_2^{j_1 + \ldots + j_4 - 1}  y_3^{j_3}  y_4^{j_4} = 0 \, .
$$
Thus one can easily see that
$$
P_{1, 1}, \ldots, P_{1, 3} \in  \operatorname{supp}({\mathfrak s_1}(H)) \quad \mbox{iff} \quad
h_{1,0,0,0} =  h_{0,1,0,0} = 0.
$$
\noindent
{\bf D$_5$}: Now  $\tilde{Y}^{1} \subset  X \times \PP_3$ has an  A$_1$ point and an A$_3$ point  as its only singularities.
One can check that
$$
P_{1,1} := (P,(1:0:0:0)) \in \operatorname{sing}(\tilde{Y}^{1})
$$
for every choice of the polynomial $F$ in \eqref{eq-sqh}.
Since the germ of $\tilde{Y}^{1}$ in $P_{1,1}$ is given by a semiquasihomogenous equation with regard to the weights
$
\operatorname{w}(x_1) = \operatorname{w}(y_3) = \operatorname{w}(y_4) = \frac{1}{2} \mbox{ and } \operatorname{w}(y_2) = \frac{1}{4},
$
it is an A$_3$ point (see \cite[Char.~C9]{durfee}).
The other singularity of  the variety $\tilde{Y}^{1}$ is the
A$_1$ point
$
P_{1,2} := (P,(-f_{0,3,0,0}:1:0:0)).
$ \\
Finally, we blow-up the singularity $P_{1,1}$. Then,
 for every polynomial $F$ of order $> 1$,
the set  $\operatorname{sing}(\tilde{Y}^{2})$ consists of the unique
A$_1$ point $P_{2,1} := (P_{1,1},(0:1:0;0))$.

By direct computation we obtain :
\begin{eqnarray}
P_{1,1} \in  \operatorname{supp}({\mathfrak s_1}(H)) & \mbox{iff} & h_{1,0,0,0} = 0 \, ,  \label{eq-pjj} \\
P_{1,2} \in  \operatorname{supp}({\mathfrak s_1}(H)) & \mbox{iff} &  h_{0,1,0,0} - f_{0,3,0,0} \cdot
h_{1,0,0,0} = 0 \, , \label{eq-pdj} \\
P_{2,1} \in  \operatorname{supp}({\mathfrak s_2}(H)) & \mbox{iff} & h_{0,1,0,0} = 0 \nonumber \, .
\end{eqnarray}

\noindent
{\bf D$_m$, $m \geq 6$}: In this case, $\operatorname{sing}(\tilde{Y}^{1})$ consists of the points
$P_{1,1}, P_{1,2}$ as for $m=5$. The point $P_{1,1}$ is a D$_{m-2}$ singularity by  \cite[Char.~C9]{durfee}, so $P_{1,2}$
in an  A$_1$ point. We arrive at the conditions \eqref{eq-pjj}, \eqref{eq-pdj}. \\
One can check that,
for $j \geq 2$,
the configuration of singularities of  $\tilde{Y}^{j}$ is independent of
the polynomial $F$. Proceeding by induction one shows that \eqref{eq-inf-eq} holds.

The proof of \eqref{eq-inf-eq} for A$_m$ and E$_m$ points is analogous, so we leave it to the reader.
\end{proof}

\noindent
Remark: If $\mbox{sing}(Y)$ consists of D$_4$ points,  then
the tangent cone $\mbox{C}_{P}Y$ of $Y$ in $P \in \sing(Y)$ consists of two $3$-planes
that meet along a $2$-plane  $\Pi_{P}$. In this case, Lemma~\ref{lem-defekt} reads
\begin{equation} \label{defectcusp}
\delta_{Y} = \dim \{ H \in V_Y
 \, : \,
\Pi_P \subset  \mbox{T}_{P}H \mbox{ for every } P \in  \sing(Y)\}   - h^0({\mathcal O}_{\pp}(2Y + K_{\pp})) + 4 \cdot  \nu,
\end{equation}
where  $\mbox{T}_{P}H$ stands for the Zariski tangent space and $\# \sing(Y) = \nu$.

\VVsect
\section{Hodge numbers of big resolutions}
\label{sec:brts}

Let  $Y \subset X$ be a three-dimensional hypersurface such that all its singularities are  A-D-E and
$\mbox{sing}(\pp) \cap Y = \emptyset$.
Let a$_m$ (resp. d$_m$, resp. e$_m$) stand for the number of singularities of $Y$ of the type A$_m$ (resp.
D$_m$, resp. E$_m$).
In this case (see \eqref{tab-infp}), the number
of singular and infinitely near singular points of the variety $Y$
is the sum
\begin{equation} \label{eq-miade}
\sum_{m \geq 1} a_m \cdot \lceil m/2 \rceil +  \sum_{m \geq 4} 2 \cdot d_m \cdot \lfloor m/2 \rfloor + 4 \cdot e_6 +
7 \cdot e_7 +  8 \cdot e_8 \, .
\end{equation}
For every $P  \in \mbox{sing}(Y)$  we choose  local  (analytic)
coordinates $x_{1,P}, x_{2,P}, x_{3,P}, x_{4,P}$ such that the
equation of the germ of $Y$ at $P$ in those coordinates is of the
form \eqref{eq-sqh} and put (see Lemma~\ref{lem-defekt})
\begin{eqnarray*}
{\mathfrak V}_{Y} & := &
\{ H \in H^0({\mathcal O}_{\pp}(2Y + K_{\pp}) \otimes  \mathcal J_{\operatorname{sing}(Y)})
\mbox{ such that the condition {\rm [II]} is fulfilled} \} \, . \nonumber
\end{eqnarray*}
Then, by Lemma~\ref{lem-defekt}, the defect of $Y$ can be expressed as
\begin{equation} \label{eq-deff}
\delta_Y :=   \dim({\mathfrak V}_{Y}) -  (h^0({\mathcal O}_{\pp}(2Y + K_{\pp})) - \mu_Y) .
\end{equation}
We have  the following
generalization of  \cite[Thm~1]{Cynk2}
for hypersurfaces with A-D-E singularities.
\begin{thm} \label{ade}
Let $\pp$ be a projective normal Cohen-Macaulay fourfold and let $Y \subset X$ be a hypersurface
with A-D-E singularities
such that  $\mbox{sing}(\pp) \cap Y = \emptyset$.
Let $\tilde{\pi} : \tilde Y \ra Y$ be the big  resolution of $Y$.
Assume that the following conditions are satisfied:
\begin{description}
\item [{[A1]}] $H^{i}({\mathcal O}_{X}(-Y))= 0$ for $i \leq 3$  and   $H^{j}({\mathcal O}_{X}(-2Y))= 0$
for $j \leq 2$,
\item [{[A2]}] $H^{2}(\OOmega^{1}_{X})=0$,
\item [{[A3]}] $H^{i}(\OOmega^{1}_{X}\otimes \mathcal {\mathcal O}_{X}(-Y))=0, \mbox{ for } i = 1, 2, 3,$
\end{description}
then
\begin{equation} \label{hjj}
  h^{1,1}(\tilde Y)=h^{1}(\OOmega^{1}_{X})+\mu_Y +\delta_Y + h^{3}({\mathcal O}_{X}(-2Y)) \, . \\
\end{equation}
Moreover, if $h^2(\op) = 0$, then
\begin{align} \label{hjd}
  &h^{1,2}(\tilde Y)=
  h^{0}(\mmm)+h^{4}(\OOmega^{1}_{X})-h^{0}({\mathcal O}_{\pp}(Y + K_{\pp}))-
  \\&\rule{1.4cm}{0cm}-h^{3}(\OOmega^{1}_{X})-
    h^{4}(\OOmega^{1}_{\pp}\otimes{\mathcal O}_{\pp}(-Y))-\mu_Y+\delta_Y, \nonumber
\end{align}
where $\mu_Y$ is the number of singularities and infinitely near singularities of $Y$ and $\delta_Y$ is the defect of $Y$.
\end{thm}

\begin{proof} We maintain the notation of Sect.~\ref{sec:brvsr}. In particular,
we consider the blow-ups  $ \sigma_{j}:  {\tilde X}^{j} \ra
{\tilde X}^{j-1}$, where $j = 1, \ldots, n$, and the big
resolution $\tilde{\pi} =  \sigma_{n} \circ \ldots \circ
\sigma_{1}$. We put $\mm=\mo_{X}(Y)$ and proceed by induction on
$n$.

\noindent {\sl n=0:} Here $\tilde{Y} = Y, \tilde{X} = X$. By
\eqref{omijaosobliwosci},  we have the   isomorphism $\mathcal
N^{\vee}_{Y|\pp}\cong \mox\otimes \mm^{-1}$ and the conormal exact
sequence
\begin{equation}  \label{ciagdoklcn}
0\lra  \mox\otimes \mm^{-1}
\lra\OOmega^{1}_{\pp}
\otimes\mo_{Y}
\lra\Omega^{1}_{Y}\lra0 \, ,
\end{equation}
so we can almost verbatim follow the proof of \cite[Thm~1]{Cynk2}:

Lemma~\ref{hiomegatensoy} and the assumptions [A2], [A3] imply
\begin{equation} \label{eq-vhn}
h^{2}(\OOmega^{1}_{X} \otimes \mo_{Y}) =  0.
\end{equation}
By Lemma~\ref{formulanadefekt}.c and \eqref{eq-vhn} the cohomology sequence associated to the conormal
sequence~\eqref{ciagdoklcn} splits and we obtain
the short exact sequence
\[0\lra H^{1}(\OOmega^{1}_{\pp} \otimes \mo_{Y})  \lra
H^{1}(\Omega^{1}_{Y})\lra H^{2}(\mox\otimes\mm^{-1})\lra0 \, .\]
Use Lemma~\ref{hiomegatensoy} to compute  $h^{1}(\OOmega^{1}_{\pp} \otimes \mo_{Y})$.
The  first formula results immediately from \eqref{eq-l1-2}.

Assume $h^2(\op) = 0$. Lemma~\ref{formulanadefekt}.b implies
$
h^{3}(\Omega^{1}_{Y}) = h^2(\mo_Y) = 0.
$
Thus, from \eqref{ciagdoklcn} and \eqref{eq-vhn},  we have the exact sequence
\[0\lra H^{2}(\Omega^{1}_{Y})\lra H^{3}(\mox\otimes\mm^{-1})\lra
H^{3}(\OOmega^{1}_{\pp}\otimes\mo_{Y})\lra0 \, . \]
Lemma~\ref{hiomegatensoy} and \eqref{eq-l1-2pp} give  the second formula for $n = 0$ .

\noindent
{\sl (n-1) $\leadsto$ n:} By Lemma~\ref{formulanadefekt}.a (resp. Lemma~\ref{hiomegatenslminone}.b)
the pair $({\tilde X}^{1}, {\tilde Y}^{1})$ satisfies [A1] (resp.~[A3]). Lemma~\ref{hiomegatenslminone}.a
yields that  ${\tilde X}^{1}$ fullfils [A2]. \\ We put
$
\ml := \mo_{{\tilde X}^{1}}({\tilde Y}^{1}) \mbox{ and } \nu : = \# \sing(Y) \, .
$
Since $\mu_Y = \mu_{{\tilde Y}^{1}} + \nu$, Lemma~\ref{hiomegatenslminone}.a implies the equality
\begin{equation} \label{eq-prbs}
h^{1}(\OOmega^{1}_{{\tilde X}^{1}})+\mu_{{\tilde Y}^{1}}  = h^{1}(\OOmega^{1}_{X})+\mu_Y \, .
\end{equation}
From the inductive hypothesis we obtain
\begin{eqnarray*}
h^{1,1}(\tilde Y) & = & h^{1}(\OOmega^{1}_{{\tilde X}^{1}})+\mu_{{\tilde Y}^{1}} +\delta_{{\tilde Y}^{1}}
+ h^{3}(\ml^{-2})  \\
  & \stackrel{\eqref{eq-prbs}}{=}  &
 h^{1}(\OOmega^{1}_{X})+\mu_Y  +\delta_{{\tilde Y}^{1}}
+ h^{3}(\ml^{-2}) \\
 & \stackrel{\operatorname{Lemma}~\ref{formulanadefekt}.\operatorname{e}}{=}  &
 h^{1}(\OOmega^{1}_{X})+\mu_Y  +\delta_{Y}
+ h^{3}(\mm^{-2}) \, ,
\end{eqnarray*}
which completes the proof of the formula for $h^{1,1}(\tilde Y)$.

If $h^{2}(\op) = 0$, then we have $h^{2}({\mathcal O}_{{\tilde X}^{1}}) = 0$. Serre duality and \eqref{vanl1} yield
\begin{equation} \label{eq-hjd-p}
h^{0}(\ml \otimes
    {\mathcal O}_{{\tilde X}^{1}}(K_{{\tilde X}^{1}})) =h^{0}(\mm \otimes
    {\mathcal O}_{\pp}(K_{\pp})) \, .
\end{equation}
Moreover, the following equality holds
\begin{equation} \label{eq-hjd-d}
h^{0}(\ml^2 \otimes {\mathcal O}_{{\tilde X}^{1}}(K_{{\tilde X}^{1}}))
 \stackrel{\eqref{eq-1-pico}}{=}  \dim V_{Y}
 \stackrel{\eqref{ry}}{=}
h^{0}(\mm^2 \otimes {\mathcal O}_{\pp}(K_{\pp})   ) - \nu  + r_Y \, .
\end{equation}
By the inductive hypothesis we have
\begin{eqnarray*}
h^{1,2}(\tilde Y) & = &
  h^{0}(\ml^2 \otimes {\mathcal O}_{{\tilde X}^{1}}(K_{{\tilde X}^{1}}))+
h^{4}\OOmega^{1}_{{\tilde X}^{1}}-h^{0}(\ml \otimes
    {\mathcal O}_{{\tilde X}^{1}}(K_{{\tilde X}^{1}}))-\\
& &-h^{3}\OOmega^{1}_{{\tilde X}^{1}}-
    h^{4}(\OOmega^{1}_{{\tilde X}^{1}}\otimes \ml^{-1})-\mu_{{\tilde Y}^{1}}+\delta_{{\tilde Y}^{1}} \\
&  \stackrel{\operatorname{Lemma}~\ref{hiomegatenslminone}}{=} &   h^{0}(\ml^2 \otimes {\mathcal O}_{{\tilde X}^{1}}(K_{{\tilde X}^{1}}))+
h^{4}\OOmega^{1}_{{X}}-h^{0}(\ml \otimes
    {\mathcal O}_{{\tilde X}^{1}}(K_{{\tilde X}^{1}}))-\\
& &-h^{3}\OOmega^{1}_{{X}}-
    h^{4}(\OOmega^{1}_{{X}} \otimes \mm^{-1})-\mu_{{\tilde Y}^{1}}+\delta_{{\tilde Y}^{1}}.
\end{eqnarray*}
Use \eqref{eq-hjd-p} and  \eqref{eq-hjd-d} to complete the proof. \end{proof}

Observe that, as in the nodal case (see \cite{Cynk2}),
using Lemma~\ref{formulanadefekt}.b one can easily compute the other Hodge numbers of ${\tilde Y}$:
\begin{align}
  &h^{1}({\mathcal O}_{\tilde Y}) = h^{1}({\mathcal O}_{X}), \quad  h^{2}({\mathcal O}_{\tilde Y}) = h^{2}({\mathcal O}_{X})\label{abc} \\
  &h^{3}({\mathcal O}_{\tilde Y}) = h^{3}({\mathcal O}_{X}) +  h^{4}({\mathcal O}_{X}(-Y)) -
  h^{4}({\mathcal O}_{X}) \, .
\end{align}

\noindent Moreover, for a smooth $X$ and an ample $Y$ with A-D-E
singularities, \cite[Cor.~6.4]{esnview} reduces the assumptions
[A1], [A2], [A3]  to the vanishing $ h^{2}(\Omega^{1}_{X}) =
h^{3}(\Omega^{1}_{X}(-Y)) = 0 \, . $

The conditions [A1], [A2], [A3] are satisfied for a simplicial, toric variety $X$
and an ample hypersurface $Y$ (for an exposition of the theory of toric varieties see \cite{oda}).

\begin{cor} \label{projective}
Let $X$ be a complete simplicial toric fourfold and let $Y \subset X$ be a hypersurface with A-D-E singularities
 such that $\sing(X) \cap Y = \emptyset$. If $\mo_{X}(Y)$ is ample, then
\begin{eqnarray*}
h^{1,1}(\tilde Y) & = & h^{1}(\OOmega^{1}_{X}) + \mu_Y + \delta_Y \, , \\
h^{1,2}(\tilde Y) & = &  h^{0}(\mmm) - h^{0}({\mathcal O}_{\pp}(Y + K_{\pp}))
    - h^{4}(\OOmega^{1}_{\pp}(-Y))-\mu_Y+\delta_Y,
\end{eqnarray*}
where $\mu_Y$ (resp. $\delta_Y$) is given by the formula \eqref{eq-miade} (resp. \eqref{eq-deff}).
\end{cor}
\begin{proof}
According to \cite[Cor.~3.9]{oda}, the variety $X$ is Cohen-Macaulay.
From \cite[Bott's Vanishing,~p.~130]{oda} (see also \cite[\S.~3]{bat93}) and Serre duality,
we obtain
$$
h^{i}({\mathcal O}_{X}(-Y)) = h^{i}({\mathcal O}_{X}(-2Y)) =   0 \mbox{ for } i \leq 3.
$$
By the same argument
$
h^{i}(\OOmega^{1}_{X}(-Y))=0, \mbox{ for } i = 1, 2, 3.
$
Now \cite[Cor.~2.8]{oda} implies the equality
$
 h^{i}(\mo_{X}) = 0  \mbox{ for } i \geq 1.
$
Finally,  \cite[Thm.~3.11]{oda} yields the vanishing
$
h^{i}(\OOmega^{1}_{X}) = 0  \mbox{ for } i \neq 1.
$
\end{proof}
\noindent
Remark: Let $(N, \Delta)$ be a fan and let  $^{\#}\hspace*{-0.5ex}\Delta(j)$ stand for the number
of $j$-dimensional cones in $\Delta$.
For $X = T_{N}\operatorname{emb}(\Delta)$, by \cite[Thm.~3.11]{oda},
we have the equality
$$
h^{1}(\OOmega^{1}_{X}) = ^{\#}\hspace*{-0.6ex}\Delta(1) - 4  \, .
$$
For a discussion how to compute the value of
$h^{4}(\OOmega^{1}_{\pp}(-Y))$ see  \cite[Thm~2.14]{materov}.

\vspace*{0.5ex}

\begin{exmp}
Let $Y \subset \PP_4 $ be  a degree-$d$ hypersurface  with A-D-E singularities, where $d \geq 3$.
Then $\mathfrak{V}_Y$ is the space of degree-$(2d - 5)$ polynomials that vanish along $\sing(Y)$ and
satisfy the condition
\eqref{eq-cond-am} (resp. \eqref{eq-cond-dm}, \eqref{eq-cond-em})
in every A$_m$ (resp. D$_m$, E$_m$) point of $Y$. Thm~\ref{ade} implies the equalities
\begin{eqnarray*}
h^{1,1}(\tilde Y)  =  1 + 2 \cdot \mu_Y +  \dim(\mathfrak{V}_{Y})  - {2d - 1 \choose 4} \, , & &
h^{1,2}(\tilde Y)  =   \dim(\mathfrak{V}_{Y}) - 5 \cdot {d  \choose 4}    \, .
\end{eqnarray*}

\noindent
In particular, let
$S_5(y_0, \ldots, y_3) = 0$ be a quintic in $\PP_3$ with
ordinary double points $P_1, \ldots, P_{\nu}$ as its only singularities and let
$L(y_0, \ldots, y_3) = 0$ be a plane that meets $S_5$ transversally.
One can check that the only singularities of the  threefold quintic
$$
Y_5 \, : \, S_5 + y_4^4 \cdot L = 0
$$
are the A$_3$ points   $P_j$, $j = 1, \ldots, \nu$.
Observe that, for a given equation of a space quintic and a plane, the latter condition can be checked
with help of Gr\"obner bases (see Remark~\ref{normform}).

Let us fix a basis $g_1, \ldots, g_{126}$ of
$H^0({\mathcal O}_{\PP_4}(5))$ and define the matrix
\[
\mbox{M}_{5} := \left[
\begin{array}{cccccc}
g_1(P_1) & \cdots &   g_{1}(P_{\nu}) & (\partial g_1/ \partial y_4)(P_{1}) & \cdots &
(\partial g_{1}/ \partial y_4)(P_{\nu}) \\
     \vdots   &    & \vdots &   \vdots       &   &  \vdots           \\
   g_{126}(P_1)  &  \cdots     & g_{126}(P_{\nu}) & (\partial g_1/ \partial y_4)(P_{1}) &   \cdots     &  \partial g_{126}/ \partial y_4(P_{\nu})
\end{array}
\right] \, .
\]
By Cor.~\ref{projective} the Hodge numbers of the big resolution of $Y_5$ are
\begin{eqnarray*}
h^{1,1}(\tilde Y_5)  =  1 + 4 \cdot \nu  - \operatorname{rank}(M_5) \quad \mbox{and} \quad
h^{1,2}(\tilde Y_5)  =  101 - \operatorname{rank}(M_5)    \, .
\end{eqnarray*}
\end{exmp}

\VVsect
\section{N-fold solids} \label{sect-nsolids}

\newcommand\B{B}
\newcommand\Y{Y}
\newcommand\bb{b}

In this section  we prove  various  generalizations of the formula
\cite[Cor.~2.32]{Clemens} for Hodge numbers of big resolutions of
double covers of $\PP_3$ branched along nodal hypersurfaces.

Let $B_d \subset \PP_3(\CC)$ be  a degree-$d$ surface given by the equation
$\bb_d = 0$, where $d > n$ is divisible by $n$,
and let $Y_d$ be the $n$-fold cyclic cover of  $\PP_3$ branched along  $B_d$. Then, $Y_d$ is a
 degree-$d$ hypersurface in the weighted
projective space $\PP := \PP(1, 1, 1, 1, \frac{d}{n})$. It is defined by
\begin{equation} \label{eq-hypers}
y_4^n - \bb_d(y_0, \ldots, y_3) = 0 \, ,
\end{equation}
where $\PP_3$ is embedded in $\PP$ by the map $ \PP_3 \ni (y_0:
\ldots :y_3) \mapsto  (y_0: \ldots :y_3 : 0) \in \PP(1, 1, 1, 1,
\frac{d}{n}). $
 Recall that $\PP$ has the unique singularity
$(0:0:0:0:1)$
 and $ K_{\PP} = {\mathcal O}_{\PP}(-4 -
\frac{d}{n}). $ Obviously $\Y_d$ does not contain the singular
point of $\PP$. One can easily see that $\PP$ is a simplicial
toric variety (see \cite[p.~35]{fulton}) and the hypersurface
$Y_d$ is ample.
 Thus the pair $(\PP, Y_d)$ satisfies the
assumptions [A1], [A2], [A3] (see also \cite[Thm~1.4.1]{dol} and
\cite[Thm~2.3.2]{dol}). Therefore, if the hypersurface  $Y_d$ has
A-D-E  singularities, then we can apply  Cor~\ref{projective}.

By \cite[Thm~2.3.2]{dol} and \cite[Cor~2.3.5]{dol} we have
$h^{1}(\OOmega^{1}_{\PP}) = 1$ and the equality
\[
h^4(\OOmega^{1}_{\PP}(-d)) = 4 \cdot h^{0}({\mathcal O}_{\PP}((n-1) \cdot \frac{d}{n} - 3)) +
 h^{0}({\mathcal O}_{\PP}(d - 4)) -  h^{0}({\mathcal O}_{\PP}((n-1) \cdot \frac{d}{n} - 4)) \, .
\]
Thus the second formula of Cor.~\ref{projective} reads
\begin{equation} \label{eq-hjdc}
h^{1,2}(\tilde{Y}_d) =
 \sum_{j=1}^{n-1} {d + j \cdot \frac{d}{n} - 1 \choose 3} - 4 \sum_{j=1}^{n-1}  {j \cdot \frac{d}{n} \choose 3} - \mu_{Y_d} + \delta_{Y_d} \, ,
\end{equation}
whereas \eqref{abc} yields
\begin{align}  \label{abc-s}
h^{1}({\mathcal O}_{\tilde Y_d}) = h^{2}({\mathcal O}_{\tilde Y_d}) = 0
\mbox{ and }
h^{3}({\mathcal O}_{\tilde Y_d}) = h^{0}({\mathcal O}_{\PP}((n-1) \cdot \frac{d}{n}-4)) \, .
\end{align}

Let $B_d$
 be a surface with  Du Val singularities and let
 a$_m$, $m \geq 1$ (resp. d$_m$, $m \geq 4$, resp. e$_m$, $m=6, 7, 8$)
stand for the number of singularities of $B_d$ of the type A$_m$ (resp.  D$_m$, resp. E$_m$).
We define the integer $\mu_{B_d}$ as the sum \eqref{eq-miade}. \\
For every   point $P \in \sing(B_d)$,
we choose  such local coordinates $x_{1,P}, x_{2,P}, x_{3,P}$ on  $\PP_3$  that the germ of $B_d$ at $P$
is locally given  by the equation
\begin{equation} \label{eq-norm3}
{\mathfrak n}(x_{1,P}, x_{2,P}, x_{3,P}) + F(x_{1,P}, x_{2,P}, x_{3,P}) = 0 \, ,
\end{equation}
where ${\mathfrak n}(x_{1,P}, x_{2,P}, x_{3,P})$ is the normal form from the table \eqref{tab-norm} and
$F(x_{1,P}, x_{2,P}, x_{3,P})$ is
 a polynomial of order $> 1$ with respect to the corresponding weights. \\
Let $d \geq 4$ be an even integer.
We define (see Lemma~\ref{lem-defekt})
\[
\begin{array}{lcl}
{\mathfrak V}_{B_d,2} & := &
\{ H \in H^0({\mathcal O}_{\PP_3}(3d/2 - 4) \otimes  {\mathcal I}_{\operatorname{sing}}(B_d))
\mbox{ such that the condition {\rm [II]} is fulfilled} \} , \nonumber  \\
\delta_{B_d,2} & := & \dim({\mathfrak V}_{B_d,2}) -  {3d/2  - 1 \choose 3} + \mu_{B_d} \, .
\end{array}
\]
With this notation we have:
\begin{cor} \label{double}
If $Y_d$ is a double solid  (i.e. $n=2$) branched along a surface $B_d$ with  Du Val singularities, then
 \[
  h^{1,1}({\tilde Y}_d)  =  1 + \mu_{B_d}  + \delta_{B_d,2}  \, , \mbox{ and }
 h^{1,2}({\tilde Y}_d)   =  { 3d/2 - 1 \choose 3 } - 4 { d/2 \choose 3} - \mu_{B_d}  + \delta_{B_d,2} \, .
\]
\end{cor}
\begin{proof}
If we put $x_{4,P} := y_4$ for $P \in \operatorname{sing}(B_d)$,
 then in the local coordinates $x_{1,P}, \ldots, x_{4,P}$
the hypersurface $Y_d$ is given by the equation~\eqref{eq-sqh}. \\
Observe that
every $H \in  H^0({\mathcal O}_{\PP}(2Y_d + K_{\PP}))$
that is divisible by $y_4$ belongs to ${\mathfrak V}_{Y_d}$. The latter implies the equality
$$
 h^0({\mathcal O}_{\PP}(2Y_d + K_{\PP}))
 - \dim({\mathfrak V}_{Y_d})
=  h^0({\mathcal O}_{\PP_3}(3 \cdot d/2 - 4)) - \dim({\mathfrak V}_{B_d,2})   \, ,
$$
so we get $\delta_{Y_d} = \delta_{B_d,2}$ (see Def.~\ref{defect}
and \eqref{ry}). Now the corollary results from \eqref{eq-hjdc}.
\end{proof}

Suppose that  $Y_d$ is the cyclic $n$-fold cover of $\PP_3$ branched along
 a nodal hypersurface $B_d$
(i.e. all singularities of $B_d$ are A$_1$ points), where $d > n$.
Since each singularity of $B_d$ endows $Y_d$ with an A$_{n-1}$ point, we have
$\mu_{Y_d} = \lceil (n-1)/2 \rceil \cdot a_1.$ We define
$$
\delta_{B_d,n} := \sum^{n-1}_{j = \lceil n/2 \rceil}
(h^0({\mathcal O}_{\PP_3}(d + j \cdot  \frac{d}{n} - 4) \otimes  {\mathcal I}_{\operatorname{sing}(B_d)})
- {d + j \cdot \frac{d}{n} - 1 \choose 3}) +  \lceil (n-1)/2 \rceil \cdot  a_1 \, .
$$
In this case, Thm~\ref{ade} yields
\begin{cor} \label{nfold}
If $Y_d$ is the cyclic $n$-fold cover of $\PP_3$ branched along
 a nodal hypersurface $B_d$ with $a_1$ singularities, then
 \[
  \begin{array}{lcl}
  h^{1,1}({\tilde Y}_d) & = & 1 +  \lceil (n-1)/2 \rceil \cdot a_1 +  \delta_{B_d,n} , \\
 h^{1,2}({\tilde Y}_d)  & = &
\sum_{j=1}^{n-1} {d + j d/n - 1 \choose 3} - 4 \sum_{j=1}^{n-1}  {j d/n \choose 3} -
  \lceil (n-1)/2 \rceil \cdot a_1 + \delta_{B_d,n}   \, .
  \end{array}
\]
\end{cor}
\begin{proof}
Put $x_{1,P} := y_4$. Observe that  for $j \geq  \lceil (n-1)/2 \rceil$ the inclusion
$$
y_4^j \cdot H^0({\mathcal O}_{\PP_3}((2n - 1 - j) \cdot  \frac{d}{n} - 4) \subset  {\mathfrak V}_{Y_d}
$$
holds,  whereas
for $j \leq  \lceil (n-1)/2 \rceil - 1$ we have
$$
y_4^j \cdot H^0({\mathcal O}_{\PP_3}((2n - 1 - j) \cdot  \frac{d}{n}   - 4)
\cap  {\mathfrak V}_{Y_d}
= y_4^j \cdot H^0({\mathcal O}_{\PP_3}((2n - 1 - j) \cdot \frac{d}{n} - 4) \otimes  {\mathcal I}_{\operatorname{sing}(B_d)}) \, .
$$
This implies the equality  $\delta_{Y_d} = \delta_{B_d,n}$.  Now \eqref{eq-hjdc} completes the proof.
\end{proof}

Finally, assume that all singularities of $B_d$ are ordinary cusps (i.e. A$_2$ points).
Recall that for $P \in \operatorname{sing}(B_d)$,
the tangent cone $\mbox{C}_{P}B_d$ consists of two planes meeting along a line.
We denote this line by  $L_P$. \\
Let $Y_d$ be the triple cyclic cover of  $\PP_3$ branched along the hypersurface
$B_d$.
Every singular point of the surface $B_d$ endows the threefold $Y_d$ with a singularity of the type D$_4$
(see  \cite[Char.~C.~9]{durfee}).
By \eqref{tab-infp} we have
 $\mu_{Y_d} = 4 \cdot a_2$. We define
$$
{\mathfrak V}_{B_d,3} := \{ H \in  H^0({\mathcal O}_{\PP_3}(5d/3 - 4)  \otimes {\mathcal I}_{\operatorname{sing}(B_d)} ) \, : \,
L_P \subset \mbox{T}_{P}H \mbox{ for every } P \in \mbox{sing}(B_d) \} \, ,  \\
$$
where $\mbox{T}_{P}H \mbox{ stands for the Zariski tangent space}$, and
$$
  \delta_{B_d,3} := \dim({\mathfrak V}_{B_d,3}) - {5d/3 - 1 \choose 3} +
h^0({\mathcal O}_{\PP_3}(4d/3 - 4)  \otimes {\mathcal I}_{\operatorname{sing}(B_d)} )
 - {4d/3 - 1 \choose 3}
 + 4 a_2.
$$
In this case, Thm~\ref{ade} implies
\begin{cor} \label{triple}
Let $Y_d$ be the triple cover of  $\PP_3$ branched along a hypersurface
$B_d$, where $d \geq 6$, with $a_2$ ordinary cusps as its only singularities. Then
 \[
 \begin{array}{lcl}
  h^{1,1}({\tilde Y}_d)  & = &  1 + 4 \cdot a_2 +  \delta_{B_d,3}, \\
 h^{1,2}({\tilde Y}_d)   & = &
{4 d/3 - 1 \choose 3} + {5d/3  - 1 \choose 3} -
4 \cdot ({2d/3  \choose 3} + { d/3 \choose 3}) - 4 a_2 +  \delta_{B_d,3}  \, .
\end{array}
\]
\end{cor}
\begin{proof}
We define $x_{4,P} := y_4$ for every $P \in \sing(Y_d)$.
We are to compute the defect $\delta_{Y_d}$.

By \eqref{defectcusp} the condition \eqref{eq-cond-dm}
 reads
\begin{equation} \label{deftripp}
H(P) =  \frac{\partial H}{ \partial x_{1,P}}(P) =  \frac{\partial H}{ \partial x_{4,P}}(P) = 0 \, .
\end{equation}
Therefore, we have the equality
$$
H^0({\mathcal O}_{\PP_3}(5d/3 - 4)) \cap   {\mathfrak V}_{Y_d} =   {\mathfrak V}_{B_d,3} \, .
$$
Observe that a hypersurface
$H \in  y_4 \cdot H^0({\mathcal O}_{\PP_3}(4d/3 - 4))$
satisfies
 \eqref{deftripp}
for every $P \in \operatorname{sing}(Y_d)$
iff we have
$H \in  y_4 \cdot H^0({\mathcal O}_{\PP_3}(4d/3 - 4)  \otimes {\mathcal J}_{\operatorname{sing}(B_d)} ).$ \\
Finally, one can easily see that  every  $H \in  H^0({\mathcal
O}_{\PP}(5d/3 - 4))$ that is divisible by  $y_4^2$ fulfills the
condition  \eqref{deftripp}.
 We obtain the equality $\delta_{Y_d}
= \delta_{B_d,3}$ and the proof is complete.
\end{proof}

We end this section with an example which shows
that working on a singular ambient variety $X$
is more efficient than dealing with a desingularization of $X$.
\newcommand\tp{\tilde{\mathbb P}}

\begin{exmp} \label{nagladkimnieidzie}

Let $Y_6 \subset \PP := {\mathbb P}(1,1,1,1,2)$ be the triple
cover of $\PP_3$ branched along a sextic surface and let $\tp$ be
the blow-up of ${\mathbb P}(1,1,1,1,2)$ in the singular point $(0:
\ldots : 1)$. Since $Y_6$ does not pass through the center of the
blow-up, its proper  transform is no longer ample.
 By abuse of
notation we use  $Y_6$ to denote the proper transform in question.
We claim that the pair $(\tp, Y_6)$ does not satisfy the
assumption [A3] of Thm.~\ref{ade}, i.e.
\begin{equation} \label{nonvanish}
h^{1}(\Omega^{1}_{\tp}(-Y_6)) = 1 \, .
\end{equation}
One can show that  $\tp = \PP({\mathcal E})$ with ${\mathcal E} := {\mathcal O}_{\PP_3} \oplus {\mathcal O}_{\PP_3}(2)$
(we maintain the notation of  \cite[Ex.~III.8.4]{hart}). Then, we have
$
K_{\tp} =(\pi^{*}{\mathcal O}_{\PP_3}(-2))(-2) \, ,
$
where $\pi : \tp \ra \PP_3$ stands for  the bundle projection
and
$\Pic( \tp ) = \ZZ^2$. From \cite[Def., p.~429]{hart}, we get
${\mathcal O}_{\tp}(1)^2 = 2\pi^{*}{\mathcal O}_{\PP_3}(1).{\mathcal O}_{\tp}(1)$, which yields
\begin{equation} \label{eq-ch1}
\pi^{*}{\mathcal O}_{\PP_3}(1)^{4} = 0 \quad \mbox{ and } \quad
   {\mathcal O}_{\tp}(1)^k.\pi^{*}{\mathcal O}_{\PP_3}(1)^{4-k} = 2^{k-1}
\mbox{ for } k = 1, \dots, 4.
\end{equation}
By  studying the canonical quotient map $\PP_4  \ra \PP$
(see \cite[App.~B]{Dimca2})
we obtain the equality
\begin{equation} \label{eq-ch2}
Y_6^k.(\pi^{*}{\mathcal O}_{\PP_3}(1))^{4-k} = 3^k \cdot 2^{k-1} \mbox{ for } k = 1, \dots, 4.
\end{equation}
We claim that $Y_6 = 3 {\mathcal O}_{\tp}(1)$. Indeed, let
$Y_6 = a \, \pi^{*}{\mathcal O}_{\PP_3}(1)  +  b \, {\mathcal O}_{\tp}(1)$.
Since $Y_6. \pi^{*}{\mathcal O}_{\PP_3}(1)^3 = 3,$
the equality \eqref{eq-ch1} with $k = 1$ implies $b = 3$. Now,  \eqref{eq-ch1} with $k = 4$ and  \eqref{eq-ch2} give $a = 0$.

By \cite[Ex.~III.8.4.c]{hart} we have
\[
\pi_{*}(\pi^{*}{\mathcal O}_{\PP_3}(-2)(1)) =  {\mathcal O}_{\PP_3}(-2) \oplus   {\mathcal O}_{\PP_3} \mbox{ and }
\operatorname{R}^{1}\pi_{*}({\mathcal O}_{\PP_3}(-2)(1)) = 0 \, ,
\]
so the Leray spectral sequence implies
$
h^i(\tilde{\PP}, \pi^{*}{\mathcal O}_{\PP_3}(-2)(1)) =
h^i(\PP_3,  {\mathcal O}_{\PP_3}(-2) \oplus   {\mathcal O}_{\PP_3}), \mbox{ where } i \geq 0 .
$ \\
Serre duality yields the equalities
\begin{equation} \label{eq-ch0}
h^4({\mathcal O}_{\tilde{\PP}}(-3)) = 1 \mbox{ and }  h^j({\mathcal O}_{\tilde{\PP}}(-3)) = 0 \mbox{ for } j \leq 3 \, .
\end{equation}
In similar way we show that
\begin{equation} \label{eq-ch3}
h^{1}(\pi^{*}{\mathcal E})(-4)  =  1,   \quad
h^{4}(\pi^{*}{\mathcal E})(-4)  =  12, \quad
h^{j}(\pi^{*}{\mathcal E})(-4)  =  0 \quad  \mbox{for } j = 0, 2, 3.
\end{equation}
Thus the exact sequence  \cite[Ex.~III.8.4.b]{hart} tensored with  ${\mathcal O}_{\tilde{\PP}}(-3)$:
$$
0 \lra \Omega_{\tilde{\PP}/{\PP_3}}(-3) \lra (\pi^{*}{\mathcal E})(-4)  \lra {\mathcal O}_{\tilde \PP}(-3)
\lra 0 \,
$$
gives the equalities
\begin{equation} \label{koh1}
h^{1}(\Omega_{\tilde{\PP}/{\PP_3}}(-3)) =  1,  \quad
 h^{4}(\Omega_{\tilde{\PP}/{\PP_3}}(-3)) =  11 \quad \mbox{and} \quad
h^{j}( \Omega_{\tilde{\PP}/{\PP_3}}(-3)) = 0 \mbox{ for } j = 0, 2, 3,
\end{equation}
In order to compute $h^{j}(\pi^{*} \Omega_{\PP_3}(-3))$, we consider  the pull-back of the Euler sequence under the map $\pi$
and tensor it with ${\mathcal O}_{\tilde{\PP}}(-3)$:
\begin{equation} \label{eq-ch4}
0 \lra (\pi^{*} \Omega_{\PP_3})(-3) \lra \oplus_{1}^{4} (\pi^{*} {\mathcal O}_{\PP_3}(-1))(-3)
\lra {\mathcal O}_{\tilde{\PP}}(-3) \lra 0 \, .
\end{equation}
We use Serre duality, \cite[Ex.~III.8.4.a]{hart} and the Leray spectral sequence to show that
\begin{equation*} \label{eq-ch5}
 h^j( \pi^{*}({\mathcal O}_{\tilde{\PP}}(-1))(-3)) = 0 \mbox{ for } j \leq 3,  \mbox{ and }
h^4(\pi^{*}({\mathcal O}_{\tilde{\PP}}(-1))(-3)) = 4.
\end{equation*}
The latter, combined with \eqref{eq-ch0} and \eqref{eq-ch4}, yields
\begin{equation} \label{eq-ch6}
h^4( \pi^{*} \Omega_{\PP_3}(-3)   ) = 15 \mbox{ and }
  h^j( ( \pi^{*}({\mathcal O}_{\tilde{\PP}}(-1))(-3))  ) = 0 \mbox{ for } j \leq 3 \, .
\end{equation}
Finally, we tensor  the exact sequence
\[
0 \lra \pi^{*} \Omega_{\PP_3}  \lra \Omega_{\tilde{\PP}} \lra  \Omega_{\tilde{\PP}/{\PP_3}}
\lra 0
\]
with ${\mathcal O}_{\tilde{\PP}}(-3)$, and apply \eqref{koh1}, \eqref{eq-ch6}  to see that \eqref{nonvanish} holds.
\end{exmp}

\VVsect
\section{Small resolutions versus big resolutions} \label{sect-svb}

Let us assume that
$\mbox{sing}(Y) = \{ P_1, \ldots, P_{\nu} \}$ consists of Gorenstein singularities.
Suppose that the threefold $Y$ has a {\sl small resolution}
$\hat{\pi} : \hat{Y}  \rightarrow Y$, i.e. $\hat{\pi}$ is a proper holomorphic map such that
 $\hat{Y}$ is smooth ,
$\hat{\pi}|_{\hat{Y} \setminus \hat{\pi}^{-1}(\operatorname{sing}(Y))}$ is an isomorphism onto the image and
the exceptional set
$$
\hat{E} :=  \hat{\pi}^{-1}(\operatorname{sing}(Y))
$$
is a curve.
By \cite[Thm~1.3]{friedman} and \cite[Thm~1.5]{friedman} (see also \cite{reid2})
the exceptional set $ \hat{\pi}^{-1}(P_l)$, where $l = 1, \ldots, \nu$, consists of
smooth rational curves meeting transversally in seven possible configurations.
In particular,  the fibers of $\hat{\pi}$ are connected. \\
It should be pointed out that a small resolution does
not have to be  K\"ahler (\cite{Werner}).
However, since  the algebraic dimension of $\hat{Y}$ is maximal, i.e. it equals three,
by \cite[Appendix~B, Thm~4.2]{hart}) if a small resolution is K\"ahler, then  it is projective.  \\
Let $\tilde{Y}$ be a smooth projective variety and let
$\tilde{\pi} : \tilde{Y}  \rightarrow Y$ be a projective morphism such that
$\tilde{\pi}|_{\tilde{Y} \setminus \tilde{\pi}^{-1}(\operatorname{sing}(Y))}$
is an isomorphism onto the image. Moreover, we assume that the exceptional set
$\tilde{E} := \mbox{Ex}(\tilde{\pi})$ is a divisor in $\tilde{Y}$.
 Observe that the big resolution of a hypersurface with  A-D-E singularities
(see Sect.~\ref{sec:brvsr}) satisfies the above conditions.

In this section we compare  Hodge numbers of the manifolds $\hat{Y}$,  $\tilde{Y}$
in the case when the former is  K\"ahler and certain cohomology groups vanish.

By \cite[Thm~1.5]{friedman} and Mayer-Vietoris we have
\begin{equation} \label{kohp1}
H^{j}(\hat{E}, \CC) =  0 \mbox{ for } j \neq 0, 2 \, .
\end{equation}
We consider the Leray spectral sequence (\cite[$\S$ 4.17]{godem} or \cite[Thm.~12.13]{mccleary}) for the constant sheaf
$\underline{\CC}_{\hat{Y}}$ and the map $\hat{\pi}$:
$$
E^{p,q}_2 = H^{p}(Y, R^q \hat{\pi}_{*} \underline{\CC}_{\hat{Y}}) \Rightarrow
H^{p+q}( \hat{Y}, \underline{\CC}_{\hat{Y}}) \, .
$$
The fibers of $\hat{\pi}$ are connected, so $\hat{\pi}_{*} \underline{\CC}_{\hat{Y}} =  \underline{\CC}_{Y}$ and we get
$
E^{j,0}_2 \cong  H^j(Y, \CC) \, .
$ \\
Since $\hat{Y}$ and $Y$ are locally compact and $\hat{\pi}$ is proper,
\cite[Remarque~4.17.1]{godem} (see also \cite[Thm~4.11.1]{godem}) yields that, for $j > 0$,  the sheaf
$R^j \hat{\pi}_{*} \underline{\CC}_{\hat{Y}}$
is a sky-scraper sheaf concentrated in the singularities of $Y$. Furthermore, we have
\begin{equation} \label{imhat}
(R^j \hat{\pi}_{*} \underline{\CC}_{\hat{Y}})_{P_l} \cong  H^j(\hat{\pi}^{-1}(P_l), \CC)
\end{equation}
for $l = 1 , \ldots, \nu$.
Thus \eqref{kohp1} means that the only non-zero $E^{p,q}_2$-term for $q > 0$ is
$$
E^{0,2}_2 \cong H^{2}(\hat{E}, \CC) \, ,
$$
and we have the long exact sequence (see e.g. \cite[Ex.~1.D]{mccleary})
\[
\ldots \ra  H^{j+2}(\hat{Y}, \CC) \ra  E^{j,2}_2 \ra  E^{j+3,0}_2 \lra H^{(j+1)+2}(\hat{Y}, \CC)
 \lra  E^{j+1,2}_2 \ra  E^{(j+1)+3,0}_2  \ra  \ldots \, .
\]
The latter yields the equalities
\begin{equation} \label{small}
h^{j}(Y, \CC ) = h^{j}(\hat{Y}, \CC) \mbox{ for } j = 0, 1, 4, 5, 6 \, .
\end{equation}

After those preparations we can prove
\begin{prop} \label{comparing}
If $h^1(\motx) = 0$, $h^3(\tilde{E},\CC) = 0$ and $\hat{Y}$ is K\"ahler, then
\[ h^{2,2}(\tilde{Y}) = h^{2,2}(\hat{Y}) +  h^4(\tilde{E}, \CC) \,.\]
\end{prop}
\begin{proof}
The varieties $\tilde{Y}$ and $\hat{Y}$ are smooth and birationally equivalent
(recall that $\hat{Y}$ is projective), so the equalities
$
h^{3,1}(\tilde{Y}) = h^{2,0}(\tilde{Y}) = h^{2,0}(\hat{Y}) = h^{3,1}(\hat{Y})
$
hold. Therefore, it suffices to show that
\begin{equation} \label{redd}
 h^{4}(\tilde{Y}, \CC) = h^4(\hat{Y}, \CC) +  h^4(\tilde{E}, \CC) \, .
\end{equation}

We consider the Leray spectral sequence
$ E^{p,q}_2 = H^{p}(Y, R^q
\tilde{\pi}_{*} \underline{\CC}_{\tilde{Y}}) \Rightarrow  H^{p+q}(
\tilde{Y}, \underline{\CC}_{\tilde{Y}}). $ As in the proof of
\eqref{imhat} we show that, for $j > 0$,
 the sheaf
$R^j \tilde{\pi}_{*} \underline{\CC}_{\tilde{Y}}$
is a sky-scraper sheaf concentrated in the singularities of $Y$ and for all $P_l \in \mbox{sing}(Y)$ we have
\begin{equation} \label{imtilde}
(R^j \tilde{\pi}_{*} \underline{\CC}_{\tilde{Y}})_{P_l} \cong  H^j(\tilde{\pi}^{-1}(P_l), \CC) \,  .
\end{equation}
Since the fibers of $\tilde{\pi}$ are connected, we have
$\tilde{\pi}_{*} \underline{\CC}_{\tilde{Y}} =  \underline{\CC}_{Y}$ and the isomorphisms
$
E^{j,0}_2 \cong H^j(Y, \CC).
$

We want to show that
\begin{equation} \label{e40infty}
E^{0,4}_{\infty} \cong E^{0,4}_2 \cong H^4(\tilde{E}, \CC)   \, .
\end{equation}
At first we check that $E^{r, 5 - r}_2$-terms vanish for  $r \in
{\mathbb Z}$. The latter is obvious for $r \neq 5$ because all
$E^{p,q}_2$-terms vanish for $p, q \neq 0$, and $E^{0,5}_2 =
H^{5}(\tilde{E}, \CC) = 0$ (recall that $\dim( \tilde{E}) = 2$).
It remains to prove that $ E^{5,0}_2 \cong H^{5}(Y, \CC) = 0 \, .
$ In view of \eqref{small}, it suffices to show that
$H^{5}(\hat{Y}, \CC) = 0$. But
 we have the vanishing $ 0 =
h^{1,0}(\tilde{Y}) = h^{1,0}(\hat{Y}) = h^{3,2}(\hat{Y}) \, . $
Since $\hat{Y}$ is K\"ahler, we obtain the equality
$h^5(\hat{Y}, \CC) = 2 \cdot h^{3,2}(\hat{Y})= 0$. \\
In particular, we have  $E^{r, 5 - r}_l = 0$ for $l \geq 2, r \in {\mathbb Z}$, 
so the differential $ d_r : E^{0,4}_r \ra E^{r, (5 - r)}_r $ is
the zero map for $r \geq 2$. The latter yields
 \eqref{e40infty}
 because the $E^{p,q}_2$-terms vanish for $p < 0$.

Now, we claim that
\begin{equation} \label{eq-httttt}
E^{4,0}_{\infty} \cong E^{4,0}_2 \cong H^4(Y, \CC) \, .
\end{equation}
Indeed, all differentials $d_r : E^{4,0}_{r} \ra E^{4 + r, -r + 1}_{r}$, where $r \geq 2$,  are trivial because
there are no non-zero $E^{p,q}_2$-terms for $q < 0$. In order to control the maps
 $d_r : E^{4-r,r-1}_{r} \ra E^{4,0}_{r}$,
one has to observe that, for $r \geq 2$,
we have
$
E^{0,3}_{r} \cong E^{0,3}_2 \cong H^3(\tilde{E},\CC) = 0  \, .
$
Thus if $r \geq 2$, then the $E^{4-r,r-1}_{r}$-term vanishes, and we obtain \eqref{eq-httttt}.

To complete the proof of \eqref{redd} observe that
$
E^{r,4-r}_{\infty} = 0 \mbox{ for } r \neq 0, 4,
$
and apply \eqref{small}.
\end{proof}

In Examples~\ref{ex-sixa},~\ref{ex-sixb} we assume
 that $\hat{Y}$ is K\"ahler and  $\mbox{h}^1({\mathcal O}_{\tilde Y}) = 0$.
\begin{exmp} \label{ex-sixa}
Let $\mbox{sing}(Y) = \{ P_1, \ldots, P_{\nu} \}$ consist of D$_4$ points.
Assume that  the variety
$Y$ is given in a neighbourhood $U_j$
of the point $P_j$ in the  local
 coordinates
$x_{1,j}, \ldots, x_{4,j}$ by the equation
\begin{equation} \label{singularity}
x_{1,j} \cdot x_{2,j} + x_{3,j}^3 + x_{4,j}^3 = 0 \, .
\end{equation}
Let  $\varepsilon$ be a primitive root of unity of order three. We put
$
f_{j,k} :=  x_{3,j} +  \varepsilon^k \cdot x_{4,j} \, .
$
For every $P_j \in \mbox{sing}(Y)$ and  $k = 0, 1, 2$, we define
the local Weil divisors on $Y$
\begin{equation} \label{eq-6-1}
D_{j,k} \, \, : \, \, x_{1,j} = 0 , \,  \, \, \,  \, \,
f_{j,k} = 0  \, .
\end{equation}
According to \cite[$\S$~2.7]{br} (see also \cite[p.~101]{Werner}), a small resolution of $U_j$
can be obtained as the projection  from the closure of the graph of the meromorphic map
$$
U_j \ni (x_{1,j}, \ldots, x_{4,j}) \mapsto ((x_{1,j}: f_{j,0}),
(x_{1,j}:  f_{j,0} \cdot  f_{j,1})) \in \PP_1 \times \PP_1 \, .
$$
One can easily see that
$\hat{\pi}^{-1}(P_j)$ consists of the
rational curves
$
(0:1) \times \PP_1, \PP_1 \times (1:0) \,.
$

Let $\sigma_1: \tilde{Y}^{1} \ra Y$ be the blow-up of $Y$ in $\mbox{sing}(Y)$.
By direct computation   $\operatorname{sing}(\tilde{Y}^{1})$
consists of three ordinary double points
on the rational curve where the components of the exceptional divisor
$(\sigma_1)^{-1}(P_j)$ meet.
After blowing up the nodes,  we get the big resolution. Thus
the exceptional divisor $\tilde{E}_{P_j} := \tilde{\pi}^{-1}(P_j)$
consists of three quadrics and two copies of $\PP_2$ blown-up
in three points. Therefore, by Mayer-Vietoris, we have
\begin{equation} \label{eulertrip}
 h^3(\tilde{E}_{P_j}, \CC) = 0, \, \,  h^4(\tilde{E}_{P_j}, \CC) = 5 \mbox{ and }
 h^{1,1}(\tilde{Y}) = h^{1,1}(\hat{Y}) +  5 \cdot \nu \,.
\end{equation}
\end{exmp}
\begin{exmp} \label{ex-sixb}
Suppose that  $\mbox{sing}(Y)$ consists of  A$_m$ points such that  all $m$'s are odd.
Let $P_j$,  where  $j = 1, \ldots, \nu$,
be  locally given by the equation:
$$
x_{1,j}^{2 (k_j + 1)} + x_{2,j}^2 + x_{3,j}^2 + x_{4,j}^2 \, .
$$
If we blow up the germ of the surface
\begin{equation} \label{eq-6-2}
x_{1,j}^{k_j + 1} - i \cdot x_{2,j} = x_{3,j} - i \cdot x_{4,j} = 0 \, ,
\end{equation}
then the proper transform is smooth and the singular point $P_j$ is replaced with a  copy of $\PP_1$
(see \cite[Ex.~2.2]{laufer}, \cite{reid2} for more details). After performing such blow-ups in every singularity, we obtain
a  small resolution  $\hat{Y}$. \\
Recall  that the exceptional divisor $\tilde{E}_{P_j}$ of the big
resolution of an A$_{2k_j + 1}$ point consists of $(k_j+1)$ smooth
rational surfaces $\tilde{E_1}, \ldots, \tilde{E}_{k_j+1}$, where
 $\tilde{E_l}, \tilde{E}_{l+1}$ meet along a smooth rational curve
and  $\tilde{E}_{l_1} \cap \tilde{E}_{l_2} = \emptyset$ for $|l_1 - l_2| > 1$.
For $l \leq k_j$  the component $\tilde{E_l}$
is a Hirzebruch surface $\mbox{F}_2$, whereas  $\tilde{E}_{k_{j}+1}$ is a smooth  quadric.
We have
\begin{equation}   \label{eulerAk}
 h^3(\tilde{E}_{P_j}, \CC) = 0, \, \,  h^4(\tilde{E}_{P_j}, \CC) = k_j + 1 \mbox{ and }
h^{1,1}(\tilde{Y}) = h^{1,1}(\hat{Y}) +  \sum_{1}^{\nu}   (k_j + 1)  \, .
\end{equation}
\end{exmp}

\begin{remark} \label{remareczka}
{\em  a) By \cite[Satz, p.~103]{Werner} the variety $Y$ in
Ex~\ref{ex-sixa} has a small resolution $\hat{Y}$ that is K\"ahler
if and only if for every $P_j$ the local Weil divisors $D_{j,k}$
(see \eqref{eq-6-1}), where $k = 0, 1, 2$, can be prolonged to
(global) Weil divisors on $Y$ that are smooth at $P_j$. In
particular, in this case the global divisors that prolong
$D_{j,k}$, for a fixed $j$ and various $k$, have no common
component through the
 point  $P_j$. \\
b)  For even $m$,
the threefold singularities A$_{m}$  have no small resolutions
by \cite[Cor.~1.16]{reid2}. }
\end{remark}


\VVsect
\section{Triple sextics and double octics} \label{sect-ts}
Here we compute the Hodge numbers of
K\"ahler small resolutions $\hat{Y}$ of triple (resp. double) covers of $\PP_3$ branched along various sextics (resp. octics)
with A$_2$ singularities (resp. A$_{j}$, where $j \geq 3$ is odd).
 We use the same symbol to denote a hypersurface
and its defining polynomial.

Consider the  manifold $\hat{Y}_{6}$
 obtained as a small resolution of  a triple cover $Y_{6}$ of $\PP_3$
branched along a sextic $B_{6}$ with A$_2$ singularities.
Obviously the canonical class $K_{Y_6}$
is trivial (see \cite[Prop.~5.73]{Kollar}).
 Since  $\hat{Y}_{6}$
is a crepant resolution, we have  $K_{\hat{Y}_6} = 0$. Assume that
$\hat{Y}_{6}$ is K\"ahler. By \eqref{abc-s} it carries neither
global $1$-forms nor global $2$-forms. Consequently,
 $\hat{Y}_{6}$ is a Calabi-Yau manifold.

 Let $\mbox{sing}(B_6) = \{ P_1, \ldots, P_{\nu} \}$.
In order to compute the Hodge numbers of $\hat{Y}_6$,
we fix a basis $h_1, \ldots , h_{35}$ of  $H^0({\mathcal O}_{\PP_3}(4))$
and a basis $f_1, \ldots , f_{84}$ of  $H^0({\mathcal O}_{\PP_3}(6))$.
We define the matrix
\[
\mbox{M}_4 := \left[ h_i(P_j) \right]_{i=1,\ldots,35}^{j =1,\ldots, \nu} \, .
\]

\noindent
For every
$P_j \in \operatorname{sing}(B_6)$,
the Hessian $\operatorname{H}_{B_6}(P_j)$ vanishes (in $\CC^4$)
along two $3$-planes
$\Pi_{1,j}, \Pi_{2,j}$.
Their common part consists of a $2$-plane.
For  $j = 1, \ldots, \nu$ we choose
a vector
 $v_j \in \CC^4$,
 such that
$
\operatorname{span}(v_j, (P_{j,1}, \ldots, P_{j,4})) =  \Pi_{1,j} \cap \Pi_{2,j}
$, and define the matrix
\[
\mbox{M}_6 := \left[
\begin{array}{cccccc}
   f_1(P_1)   &  \cdots     &  f_1(P_{\nu})    &     f'_1(P_1).v_{1}   & \cdots &   f'_1(P_{\nu}).v_{\nu}    \\
     \vdots   &             &    \vdots &       \vdots       &        &         \vdots       \\
 f_{84}(P_1)  &  \cdots     & f_{84}(P_{\nu})&   f'_{84}(P_1).v_{1}  &   \cdots     &   f'_{84}(P_{\nu}).v_{\nu}   \\
\end{array}
\right] \, .
\]
Observe that, by \eqref{abc-s} and Ex.~\ref{ex-sixa},
the variety $\tilde{Y}_6$  fulfills the assumptions of  Prop.~\ref{comparing}.
From the latter and Cor.~\ref{triple} we  obtain
\begin{cor} \label{ts}
The Hodge numbers of the Calabi-Yau manifold  $\hat{Y}_{6}$ are given by the formulae:
 \[
  \begin{array}{lcl}
  h^{1,1}(\hat{Y}_6) & = & 1 + 3 \cdot  \nu  -  \operatorname{rank}( \mbox{M}_{4})  -  \operatorname{rank}( \mbox{M}_{6}), \\
  h^{1,2}(\hat{Y}_6)  & = & 103 -  \operatorname{rank}( \mbox{M}_{4}) -  \operatorname{rank}( \mbox{M}_{6}) \, .
   \end{array}
 \]
\end{cor}
\begin{proof}
We maintain the notation of the previous section.
 For $F \in  H^0({\mathcal O}_{\PP_3}(6))$,
the condition \eqref{deftripp} means that the point $P$  belongs to $F$
 and
the (Zariski) tangent space of the cone over $F$  at the point
$P_j \in \CC^4$ contains the $2$-plane $ \Pi_{1,j} \cap
\Pi_{2,j}$. The latter amounts to  the equalities $ F(P_j) =
F'(P_j).v_{j} = 0 \, , $
 so  $\dim({\mathfrak V}_{B_6}) = 84 -  \operatorname{rank}( \mbox{M}_{6})$ and we are in position to compute the defect
\[
 \delta_{B_d} = 4 \cdot \nu -   \operatorname{rank}( \mbox{M}_{4})  -  \operatorname{rank}( \mbox{M}_{6}) \, .
\]
Thus Cor.~\ref{triple} combined with Prop.~\ref{comparing} and
\eqref{eulertrip}
 yields the formula for  $h^{1,1}(\hat{Y}_6)$.

Finally, we claim that
\begin{equation} \label{rameau}
h^{1,2}(\hat{Y}_6) = h^{1,2}(\tilde{Y}_6) \, .
\end{equation}
Indeed, since $\operatorname{e}(\tilde{\pi}^{-1}(P_j)) =  13$
(see Ex.~\ref{ex-sixa}), we have the equality
$
 \operatorname{e}(\tilde{Y}_6) =  \operatorname{e}(\hat{Y}_6 ) + 10 \cdot \nu .
$
From
\eqref{abc-s} we obtain
$
\operatorname{e}(\tilde{Y}_6) = 2 (  h^{1,1}(\tilde{Y}_6) -  h^{1,2}(\tilde{Y}_6)).
$
Therefore,  the formula for $h^{1,1}(\hat{Y}_6)$ implies \eqref{rameau}.
\end{proof}
 \noindent
Remark: As  an immediate consequence of Cor.~\ref{ts} we obtain
the equality $ \operatorname{e}(\hat{Y}_6 ) = 6 \cdot \nu - 204. $
The  latter results also  from \cite[Lemma~3]{br}. Indeed, by
Noether's formula, we have $\operatorname{e}(B_6) = 108 - 2 \cdot
\nu.$
  Thus we obtain
$\operatorname{e}(\PP_3 \setminus B_6) = 2 \cdot \nu  - 104$, and
the equality in question follows.


\vskip2mm
Now we are in position to compute the Hodge numbers of K\"ahler small resolutions of triple sextics
 branched along the surfaces discussed in \cite{br1}, \cite{br2} and \cite{labs2}.

\noindent
{\bf Direct construction of \cite{br2}:} We choose surfaces
$S_1,...,S_k$ of degrees
$$\deg(S_1)=d_1,..., \deg(S_k)=d_k, \quad d_1+...+d_k=6, \quad
$$
a quadric $S$
and consider the sextic $B_6 \subset \PP_3$ given by
the equation
\begin{equation} \label{directconst}
S_1 \cdot ... \cdot S_k - S^3 = 0 \, .
\end{equation}
We require that

\vskip3mm
(\md$_1$) any three surfaces $S_i, S_j, S$ meet transversally,

(\md$_2$) no four surfaces $S_i, S_j, S_m, S$ meet,

(\md$_3$) the surface $B_6$ is smooth away from the cusps $P_{\nu}$ at
the intersections
$S_i \cap S_j \cap S$.

\vskip3mm
\noindent
The resulting sextic \eqref{directconst} has $2 \cdot \sum_{i \neq j} d_i \cdot d_j$ cusps and no other singularities.

\begin{lem} \label{exsofsmall}
If $Y_6$ is the triple cover of $\PP_3$  branched along the sextic $B_6$
obtained by the direct construction {\em \eqref{directconst}}, then there exists a small resolution
 $\hat{Y}_{6}$ of $Y_6$ that is K\"ahler.
\end{lem}
\begin{proof}
Let  $Y_6 \subset  \PP(1, 1, 1, 1, 2)$
be given by the equation \eqref{eq-hypers} and let $P_{\nu}$ be a cusp  in $S_i \cap S_j \cap S$.
We consider the  (global) Weil divisors
\begin{equation} \label{weildiv}
W_{l} \, : \,  S_{i} = 0 , \,  \, \, \,  \, \,
 \varepsilon^l \cdot S + y_{4} = 0 \, ,
\end{equation}
where  $l = 0, 1, 2$ and  $\varepsilon$ is a primitive root of unity of order three.
Obviously, $W_l$ prolongs the local divisor $D_{\nu,l}$ (see \eqref{eq-6-1}) and
the germ of $W_l$ in the point $P_{\nu}$ is smooth.
Thus
the assumptions of  \cite[Satz, p.~103]{Werner} (see also Remark~\ref{remareczka}.a) are satisfied.
\end{proof}

To check that the conditions (\md$_1$) - (\md$_3$) are satisfied by given surfaces $S_1, \ldots, S_k, S$ we will use
the following

\begin{remark} \label{normform}
{\em In order to show that a polynomial $g \in \KK[y_0, y_1, y_2, y_3]$
belongs  to an ideal ${\mathcal I}$,
one applies the notion of the
 remainder on division
of a polynomial  $g$
 by a  Gr\"obner basis  ${\mathcal B}$   of the ideal
${\mathcal I} \subset \KK[y_0, y_1, y_2, y_3]$
(see \cite[II.$\S$6]{clos}).
It is well-known that if the remainder
vanishes, then $g$ is an element of ${\mathcal I}$.
The former can be checked e.g. with the Maple command:
normalf(g, ${\mathcal B}$,  tdeg($y_0$, $y_1$, $y_2$, $y_3$)).
If the output is zero, then $g \in {\mathcal I}$. }
\end{remark}

\newcommand\xj{y_{0}}
\newcommand\xczt{y_{1}}
\newcommand\xdw{y_{2}}
\newcommand\xtrzy{y_{3}}
\begin{exmp}  \label{example-71}
We consider the quadric
$$
S \, : \,  \xj \cdot \xczt  - \xdw \cdot  \xtrzy = 0,
$$
the planes
\begin{equation*}
\begin{array}{lll}
& F_1 \, : \,  \xj = 0,  &  F_2 \, : \,  \xczt = 0 \, , \\
& F_3 \,  : \,  4  \cdot \xj - 2  \cdot \xdw - 2  \cdot \xtrzy + \xczt = 0,
& F_4 \,  : \,  \xj - 2  \cdot \xdw - 2  \cdot \xtrzy + 4 \cdot \xczt = 0 \, , \\
& F_5 \, : \, \xj +   \xdw + \xtrzy +  \xczt = 0 \, ,
& F_6 \, : \,  \xj -   \xdw - \xtrzy +  \xczt = 0 \, ,
\end{array}
\end{equation*}
and define the sextic
$$
B_6 \, : \, F_1 \cdot  \ldots \cdot F_6 - S^3 = 0 \, .
$$
The lines $F_i = F_j = 0$, where $i \neq j \leq 6$, meet the quadric $S$ in two points, so the condition (\md$_1$)
is satisfied.  \\
Since a
Gr\"obner basis computation with Maple (see Remark~\ref{normform}) shows that the polynomials
$y_0^{10}, y_1^{10}, y_2^{10}, y_3^{10}$ belong to the ideal generated by  $F_i, F_j, F_k, S$ with $i \neq j \neq k$,
the condition (\md$_2$) is fulfilled. \\ Finally, a similar Gr\"obner basis argument shows the polynomial
$S^{10}$ belongs to the
jacobian ideal of $B_6$. Therefore, if $\operatorname{mult}_{P}(B_6) \geq 2$,
then the inequality
 $
\operatorname{mult}_{P}(F_1 \cdot \ldots F_6) \geq 2 \,
$
holds. The latter shows that  the condition
 (\md$_3$) is satisfied.

Let $\hat{Y}_{6}$ be a K\"ahler small resolution of the triple sextic branched along $B_6$.
A Maple computation yields
$
\operatorname{rank}( \mbox{M}_{4}) = 25
\mbox{ and  } \operatorname{rank}( \mbox{M}_{6}) = 55.
$
The sextic $B_6$ has $30$ singularities, so
we obtain the equalities
\[
 h^{1,1}(\hat{Y}_6) =  11 \mbox{ and }
 h^{1,2}(\hat{Y}_6) =  23 \, .
\]

\end{exmp}

\begin{exmp} \label{example-72}
We maintain the notation of Example~\ref{example-71}.  Consider the quadric
$$
 R \, : \, \xj^2 +  \xczt^2 +  \xdw^2 + \xtrzy^2 = 0
$$
and the planes
\begin{equation*}
\begin{array}{lll}
& R_1 \, : \,  \xj + 2 \cdot \xdw + 3 \cdot \xtrzy + 4 \cdot \xczt = 0 \, ,
& R_2 \, : \,  4 \cdot \xj + 3 \cdot \xdw + 2 \cdot \xtrzy + \xczt = 0\, .
\end{array}
\end{equation*}
The Gr\"obner basis computation
shows that $y_0^{10}, \ldots,  y_3^{10}$ belong to the jacobian ideals
of the surfaces
\begin{align}
& S_{1,2} \, : \,  F_1 \cdot F_2 + 2 \cdot S = 0,
 &  S_{2,3} \, : \,  F_2 \cdot F_3 + 2 \cdot S = 0, \mbox{ \hspace*{5.3ex} }  \nonumber  \\
& S_{3,4} \, : \,   F_3 \cdot F_4 + 2 \cdot S = 0,
 &  S_{5,6} \, : \,  F_5 \cdot F_6 + 2 \cdot S = 0, \mbox{ \hspace*{5.3ex} } \nonumber \\
& S_{1,2,3} \, : \, F_1 \cdot F_2 \cdot F_3 + S \cdot R_1 = 0,
&  S_{4,5,6} \, : \,  F_4 \cdot F_5 \cdot F_6 + S \cdot R_2 = 0,   \label{direct} \\
& S_{3,4,5,6} \, : \,  F_3  \cdot F_4  \cdot F_5  \cdot F_6 + S  \cdot R = 0, & \nonumber \\
& S_{2,3,4,5,6} \, : \, F_2 \cdot F_3  \cdot F_4 \cdot F_5  \cdot F_6 + R_1 \cdot S  \cdot R = 0,  & \nonumber
\end{align}
so the latter are smooth.  We obtain the following table:

\vskip3mm
\renewcommand\arraystretch{1.2}
\begin{center}
\begin{tabular}{|l|c|c|c|c|c|c|}
\hline
Equation of $B_6$
& $\sharp(\operatorname{sing}(B_6))$ &  $\operatorname{rk}(\mbox{M}_4)$ &
  $\operatorname{rk}(\mbox{M}_6)$ &  $h^{1,1}(\hat{Y}_6)$ &
$h^{1,2}(\hat{Y}_6)$   \\ \hline
$F_1 \cdot S_{2,3,4,5,6} - S^3$  & 10    & $9$ & $19$  & 3 & 75 \\  \hline
$S_{1,2} \cdot S_{3,4,5,6} - S^3$  & 16    & $15$ & $31$  & 3 & 57  \\  \hline
$S_{1,2,3} \cdot S_{4,5,6} - S^3$  & 18    & $17$ & $35$  & 3  &  51 \\  \hline
$F_1 \cdot F_{2} \cdot S_{3,4,5,6}  - S^3$  &  18  &  $16$ & $34$  & 5  & 53 \\  \hline
$F_1 \cdot S_{2,3} \cdot S_{4,5,6} - S^3$  & 22  &  $20$ & $42$  & 5  & 41 \\  \hline
$S_{1,2} \cdot S_{3,4} \cdot S_{5,6} - S^3$  & 24  &  $22$ & $46$  & 5 & 35 \\  \hline
$F_1 \cdot F_2 \cdot F_3 \cdot S_{4,5,6} - S^3$ & 24  &  $21$ & $45$  & 7  & 37  \\  \hline
$F_1 \cdot F_2 \cdot S_{3,4} \cdot S_{5,6} - S^3$ & 26 &  23 & 49  & 7 &  31 \\ \hline
$F_1 \cdot \ldots  \cdot F_4  \cdot S_{5,6} - S^3$ & 28 &  $24$ & $52$ & 9 & 27  \\ \hline
\end{tabular}
\end{center}
\renewcommand\arraystretch{1}

\vskip4mm
\noindent
For the sextics in the first column of the table,  we check that the conditions
(\md$_1$) -- (\md$_3$) are fulfilled in the way shown in the previous example. In particular,
all  singularities  of each sextic lie on the quadric $S$ and
the surfaces (\ref{direct}) are smooth, so (\md$_3$) holds. To find
the singular points of the surfaces from the table, we use the fact that they are singularities of
the sextic $B_6$ from Example~\ref{example-71}. \\
\end{exmp}

\noindent
{\bf Residual construction of \cite{br2}:}
To construct another example we apply the residual construction of \cite{br2}.
We choose

- a residual cubic $R$;

- auxiliary planes $R_1,...,R_3$
such that the curves $R_i \cap R$ are smooth and
intersect transversally,

- cubic surfaces $S_i \, : \, R_i^3 + \lambda_i \cdot R = 0$,

- a cubic $S \, : \, R_1 \cdot ... \cdot R_3 + \lambda \cdot R = 0$,

\vskip2mm
\noindent
where $\lambda, \lambda_i \in \CC$. Then, the polynomial $S_1 \cdot S_2 \cdot S_3 - S^3$ always vanishes along the residual cubic
$R$ and we can
consider the following sextic
\begin{equation} \label{residualconstr}
B_6 \, : \,  (S_1 \cdot S_2 \cdot S_3 - S^3)/R = 0 \,  .
\end{equation}
By \cite[Sect.~1.2]{br2},
if we choose the constants $\lambda, \lambda_i$ general enough, then

\vskip2mm
(\mr$_1$) $B_6$ has no singularities along the residual sextic $R$,

(\mr$_2$) the cubics $S_i, S_j, S$ intersect transversally outside of $R$,

(\mr$_3$) the sextic $B_6$ is smooth away from the points in $S_i \cap S_j \cap S$.

\vskip2mm
\noindent
The conditions (\mr$_1$) -- (\mr$_3$) imply that the sextic $B_6$ has $27$ A$_2$ singularities.

Let $Y_6 \subset  \PP(1, 1, 1, 1, 2)$
be the triple cover of $\PP_3$ branched along the sextic \eqref{residualconstr} and let
$P$ be a cusp in $S_i \cap S_j \cap S$.
As in the
proof of Lemma~\ref{exsofsmall}, we show that the divisors \eqref{weildiv} are smooth
in $P$, so
the triple sextic
$Y_6$ has a
K\"ahler small resolution $\hat{Y}_6$.
\begin{exmp} (c.f. \cite[Sect.~2.3]{br2}) \label{ex-resid}
We take as the residual cubic the
Fermat
cubic
$$R \, : \, y_0^3+y_1^3+y_2^3+y_3^3 = 0
$$
and put
$$
S_i \, : \, y_i^3+ R = 0, \, i=1,2,3, \quad S \, : \, y_1 \cdot y_2\cdot y_3 = 0.
$$
Using Gr\"obner bases (see Remark~\ref{normform})
we show that the above defined hypersurfaces and the sextic $B_6$ given by the equation \eqref{residualconstr}
satisfy the conditions
(\mr$_1$) -- (\mr$_3$). In particular,  the surface $B_6$ is smooth away from the points
$$
(1: 0: \epsilon^{i_1}:  \epsilon^{i_2}), \quad
(1: \epsilon^{i_1}: 0:  \epsilon^{i_2}), \quad
(1: \epsilon^{i_1} :  \epsilon^{i_2}: 0) \, ,
$$
where $0 \leq i_1, i_2 \leq 2$ and $\epsilon$
is a primitive root of $(- 1/3)$ of order three.
We obtain  $\operatorname{rank}( \mbox{M}_{4}) = 24 $
and $\operatorname{rank}( \mbox{M}_{6}) = 51$, so
\[
 h^{1,1}(\hat{Y}_6) =  7 \mbox{ and }
 h^{1,2}(\hat{Y}_6) =  28 \, .
\]
\end{exmp}

\noindent
{\bf Sextic with 36 cusps (see \cite{labs2}):}
Recall that, by \cite{v}, \cite{m}, the number of A$_2$ singularities on a sextic in $\PP_3$ does not exceed $37$
and it is not known whether this bound is sharp.
As the last example of a triple sextic, we consider the cover branched
along the sextic with $36$ ordinary cusps that was
constructed in \cite[App.~A]{labs2}.
\begin{exmp} \label{labssextic2}
Consider the quadric $S \, : \, y_0 \cdot y_1 - y_2 \cdot y_3 = 0$ and the linear forms
\begin{equation*}
\begin{array}{lll}
& z_0 := y_0,  &  z_2 := y_0 + y_1 - y_2 - y_3 \\
& z_1 := y_1 ,
& z_3 := 8 y_0 + 8 y_1 - 64 y_2 - y_3 \, ,
\end{array}
\end{equation*}
that define the planes   tangent to $S$ in the points
\[
(0:1:0:0), (1:0:0:0), (1:1:1:1), (8:8:1:64) \, .
\]
Let $B_6$ be the pull-back of  $S$ under the map
$
\Omega^{3}_{3} : (z_0 : z_1 : z_2 : z_3) \lra (z_0^3 : z_1^3 : z_2^3 : z_3^3).
$
Then, $B_6$ is the sextic given by the polynomial
\begin{equation} \label{labssextic}
(z_0 \cdot z_1)^3 - (8/9 z_0^3 + 8/9 z_1^3 - 64/63 z_2^3 +  1/63 z_3^3) \cdot ( 1/9 z_0^3 + 1/9 z_1^3 + 1/63 z_2^3 -  1/63 z_3^3) \, ,
\end{equation}
with $36$ A$_2$-points
\[
(\Omega^{3}_{2})^{-1}( \{ (0:1:1:8), (1:0:1:8), (1:1:0:(-49)), (8:8:(-49):0) \} )
\]
and no other singularities.
In order to show that the triple cover $Y_6 \subset  \PP(1, 1, 1, 1, 2)$
 branched along $B_6$ has a
small resolution $\hat{Y}_6$ that is K\"ahler, we use \cite[Satz, p.~103]{Werner} (see also Remark~\ref{remareczka}.a):

Obviously $B_6$ is given by the equation of type \eqref{direct} with
the quadric $z_0 \cdot z_1$.
Let $W_l$ be the divisors defined by \eqref{weildiv}, where $l = 0, 1, 2$.
As in the proof of Lemma~\ref{exsofsmall} we show that  $W_l$  are smooth in the cusps
$(\Omega^{3}_{2})^{-1}( \{ (0:1:1:8), (1:0:1:8) \}).$ \\
  If we put
$$
S' = -3^{4/3} \cdot z_2 \cdot z_3 \, , \quad S'_1 =  7 \cdot z_0^3 - 56 \cdot z_1^3 - 8 \cdot  z_2^3 -  z_3^3 \, , \quad
S'_2 =  56 \cdot z_0^3 - 7 \cdot z_1^3 + 8 \cdot  z_2^3 +  z_3^3 \, ,
$$
then, by direct computation, the equation \eqref{labssextic} can be written as
$$
(-1/3969) \cdot (S'_1 \cdot S'_2 - (S')^3).
$$
One can easily see that the divisors
$$
W'_l   \, : \,  S'_{1} = 0 , \,  \, \, \,  \, \,
 \varepsilon^l \cdot S' + y_{4} = 0 \, , \quad l = 0, 1, 2,
$$
are smooth in the cusps
$
(\Omega^{3}_{2})^{-1}( \{ (1:1:0:(-49)), (8:8:(-49):0) \} ),$
so  there exists a K\"ahler small resolution.

In this case, the Maple computation yields
 $\operatorname{rank}( \mbox{M}_{4}) =  30 $
and $\operatorname{rank}( \mbox{M}_{6}) = 66$. Therefore, we have
\[
 h^{1,1}(\hat{Y}_6) =  13 \mbox{ and }
 h^{1,2}(\hat{Y}_6) =  7\, .
\]
\end{exmp}

\begin{remark}
{\em a) For a Calabi-Yau manifold
$\hat{Y}$
both $h^{1,1}(\hat{Y})$ and $h^{1,2}(\hat{Y})$
have a geometric interpretation:
the latter is the dimension of the space of infinitesimal deformations
of the manifold in question, whereas the former equals the rank of  $\operatorname{Pic}(\hat{Y})$. \\
Observe that
the group $\operatorname{Pic}(\hat{Y}_6)$ is free. Indeed, if we repeat the proof of \eqref{small}
for the constant sheaf  $\underline{\ZZ}_{\hat{Y}_6}$, then we obtain
 the exact sequence
\[
0 \lra  H^{2}(Y_6, \ZZ) \lra H^{2}(\hat{Y}_6, \ZZ) \lra  H^{2}(\hat{E}, \ZZ) \lra
 H^{3}(Y_6, \ZZ) \lra H^{3}(\hat{Y}_6, \ZZ) \lra 0 \, .
\]
Since $H^{2}(\hat{E}, \ZZ)$ is torsion-free,
$\operatorname{Torsion}(H^{2}(\hat{Y}_6, \ZZ))$ comes from the
group $H^{2}(Y_6, \ZZ)$. But we have
 $
\operatorname{Torsion}(H^{2}(\hat{Y}_6, \ZZ)) \cong
 \operatorname{Torsion}(H_{1}(\hat{Y}_6, \ZZ)).
$ Therefore,  if   $Y_6$ is simply-connected, then
$H^{2}(\hat{Y}_6, \ZZ)$ is torsion-free. In particular, for a
weighted projective hypersurface  $Y_6$, by
\cite[Cor.~B.21]{Dimca2}, we obtain
$$
 \operatorname{Torsion}(H^{2}(\hat{Y}_6, \ZZ)) = 0 \, .
$$

\noindent
b) Recall
that the {\sl extended code} $\mathfrak{E}_{B_6}$ of  a sextic $B_6$ with
cusps $P_1, \ldots, P_{\nu}$
is defined as the kernel of the
 $\FF_3$-linear
morphism
$$
\FF_3^{\nu+1} \to H^2(\tilde{B}_6 , \FF_3), \quad
(t_0, t_1,...,t_{\nu}) \mapsto
\OO_{\tilde{B}_6}(t_0)+ \sum_{j=1}^{\nu} t_{j}[C_{j}'-C_{j}''],
$$
where $\tilde{B_6}$ is the minimal resolution of the surface $B_6$ and $C_{j}',C_{j}''$
are the exceptional $(-2)$-curves over the cusp $P_{j} \in B_6$
(see \cite[Sect.~1.3]{br2} for more details).
By \cite[Prop.~2.1]{br2} and \cite[Thm~2.9]{br2}
the extended code $\mathfrak{E}_{B_6}$ of a sextic $B_6$
given by \eqref{directconst}
depends only on the partition $d_1, \ldots, d_k$.
Furthermore, we have
$
\dim(\mathfrak{E}_{B_6}) = k - 1 \, .
$
Thus for the sextics in Ex.~\ref{example-71},~\ref{example-72} the following equality holds
$$
h^{1,1}(\hat{Y}_6) = 2 \cdot \dim(\mathfrak{E}_{B_6}) + 1.
$$
The dimension of $\mathfrak{E}_{B_6}$ is unknown for the surfaces from Ex.~\ref{ex-resid},~\ref{labssextic2}. }
\end{remark}

\vskip4mm
Consider a small resolution of the double octic i.e.
the  manifold $\hat{Y}_{8}$  obtained as a small resolution of  the double cover $Y_{8}$ of $\PP_3$
branched along an octic $B_{8}$ with A$_m$ singularities such that $m$'s are odd.
Assume that   $\hat{Y}_{8}$ is K\"ahler. As in the case of triple sextics
$\hat{Y}_6$, one can show that the
projective manifold $\hat{Y}_{8}$ is a Calabi-Yau manifold:
 the canonical class   $K_{\hat{Y}_8}$ is trivial,
and, by \eqref{abc-s}, $\hat{Y}_8$ carries neither  global $1$-forms nor global
$2$-forms.

In order to give an explicit formula for  Hodge numbers of  $\hat{Y}_{8}$, we choose a basis
$g_1, \ldots g_{165}$ of  $H^0({\mathcal O}_{\PP_3}(8))$.
For an  A$_{2k+1}$ point $P_j \in \operatorname{sing}(B_8)$ we put
\[
\mbox{M}_{8,j} := \left[
\begin{array}{ccc}
   g_1(P_j)   &  \cdots     &  (g_1|_{L_j})^{(k)}(P_{j})      \\
     \vdots   &             &    \vdots           \\
   g_{165}(P_j)  &  \cdots     & (g_{165}|_{L_j})^{(k)}(P_{j})  \\
\end{array}
\right],
\]
where the line $L_j$ is the singular locus of the set of zeroes of the Hessian  $\operatorname{H}_{B_8}(P_j)$
and $(-)^{(k)}$ denotes the $k$-th derivative.
Let  $a_{2k+1}$ stand for the  number of A$_{2k+1}$ points of the octic B$_8$.
We define  the $(165 \times (\sum a_{2k+1} \cdot (k+1)))$-matrix
\[
M_8 := [M_{8,1}, \ldots, M_{8,\nu}] \, .
\]
Then, Cor.~\ref{double} implies
\begin{cor} \label{do}
The Hodge numbers of the Calabi-Yau manifold $\hat{Y}_8$ are given by the formulae:
 \[
   \begin{array}{lcl}
   h^{1,1}(\hat{Y}_8) & = & 1 +  \sum a_{2k+1} \cdot (k+1)  -  \operatorname{rank}( \mbox{M}_{8}) \, , \\
  h^{1,2}(\hat{Y}_8)  & = & 149 -  \operatorname{rank}( \mbox{M}_{8}) \, .
   \end{array}
 \]
\end{cor}
\begin{proof} We maintain the notation of Ex~\ref{ex-sixb}. Observe that
 Cor.~\ref{double} combined with
Prop.~\ref{comparing} and  \eqref{eulerAk} gives the first formula.
By direct computation (see Ex.~\ref{ex-sixb}) we have
$\operatorname{e}(\tilde{\pi}^{-1}(P_j)) = 4 + 2 k_j$
which yields the equality
$$
 \operatorname{e}(\hat{Y}_8) =  \operatorname{e}(\tilde{Y}_8 ) - \sum a_{2k+1} \cdot (k+1) \, .
$$
From  \eqref{abc-s} we obtain $h^{1,2}(\hat{Y}_8) = h^{1,2}(\tilde{Y}_8)$
(see the proof of Cor~\ref{ts}).
\end{proof}

\noindent
Observe that in this case, \cite[Lemma~3]{br}  implies the equality
$
\operatorname{e}(\hat{Y}_8 ) = 2 \cdot \sum a_{2k+1} \cdot (k+1) - 296 \,  .
$

By \cite{m} there are no octics in $\PP_3$ with more than $69$ A$_3$ points.
The best known example is an octic $B_8$  with $64$ such singularities (\cite[App.~A]{labs2}).
Now we apply Cor.~\ref{do} to compute Hodge numbers of a K\"ahler
small resolution $\hat{Y}_8$ of the double solid branched along
the surface $B_8$.
\begin{exmp}
We maintain the notation of Ex.~\ref{labssextic2} and consider the pull-back of the quadric $S$
under the map
$$
\Omega^{3}_{4} : (z_0 : z_1 : z_2 : z_3) \lra (z_0^4 : z_1^4 : z_2^4 : z_3^4).
$$
The resulting surface is smooth away from the $64$ A$_3$ points
\[
(\Omega^{3}_{4})^{-1}( \{ (0:1:1:8), (1:0:1:8), (1:1:0:(-49)), (8:8:(-49):0) \} ) \, .
\]
The Maple computation yields
 $\operatorname{rank}( \mbox{M}_{8}) = 122$.
Therefore, if $Y$ has a small resolution that is K\"ahler, then the Hodge numbers are
\[
 h^{1,1}(\hat{Y}_8) =  7 \mbox{ and }
 h^{1,2}(\hat{Y}_8) =  27 \, .
\]
\end{exmp}

\noindent {\bf Acknowledgement:} The author would like to thank
Prof.~W.~P.~Barth for numerous fruitful discussions. We thank
Prof.~D.~van~Straten for \cite{vS}
 and  Prof.~S.~Cynk
for inspiring discussions on the papers \cite{Cynk}, \cite{Cynk2}, \cite{Cynk3}, \cite{Cynk4}. \\
The paper contains results from the author's Habilitationsschrift. The author would like to thank the Institute of Mathematics of Erlangen-N\"urnberg University
for creating perfect research conditions.

\VVsect

\vspace*{3ex}
\noindent
S{\l}awomir Rams \\
Mathematisches Institut,
Universit\"at Erlangen-N\"urnberg,
Bismarckstra\ss e 1 1/2,
D-91054 Erlangen,
Germany \\
and \\
Institute of Mathematics,
Jagiellonian University,
ul. Reymonta 4,
30-059 Krakow,
Poland

\parindent=0cm
\end{document}